\documentclass[twocolumn]{autart}    

\usepackage{amsmath,amsfonts,bm}
\usepackage{algorithm}
\usepackage{algcompatible}
\usepackage{amssymb}
\DeclareMathOperator*{\argmin}{argmin}  
\usepackage{array}
\usepackage{color}
\usepackage[caption=false,font=normalsize,labelfont=sf,textfont=sf]{subfig}
\usepackage{textcomp}
\usepackage{stfloats}
\usepackage{url}
\usepackage{verbatim}
\usepackage{graphicx}
\usepackage{cite}
\usepackage{multirow}

\newtheorem{definition}{Definition}
\newtheorem{lemma}{Lemma}
\newtheorem{remark}{Remark}
\newtheorem{corollary}{Corollary}
\newtheorem{assumption}{Assumption}          

\begin{document}

\begin{frontmatter}

\title{Distributed Empirical Risk Minimization With Differential Privacy\thanksref{footnoteinfo}} 

\thanks[footnoteinfo]{
This work was supported by the Swedish Research Council, the Knut and Alice Wallenberg Foundation, and the Swedish Foundation for Strategic Research. A preliminary version of this work has been reported in the 9th IFAC Workshop on Networked Systems \cite{LIU202243}.
}

\author[KTH]{Changxin Liu}\ead{changxin@kth.se},    
\author[KTH]{Karl H. Johansson}\ead{kallej@kth.se},               
\author[UVic]{Yang Shi}\ead{yshi@uvic.ca}  

\address[KTH]{School of Electrical Engineering and Computer Science, KTH Royal Institute of Technology, and Digital Futures, 100 44 Stockholm, Sweden}  
\address[UVic]{Department of Mechanical Engineering, University of Victoria, 
  Victoria, B.C. V8W 3P6, Canada}             

\begin{keyword}                           
Distributed optimization; empirical risk minimization; differential privacy; dual averaging.
\end{keyword}                             

\begin{abstract}                          
This work studies the distributed empirical risk minimization (ERM) problem under differential privacy (DP) constraint. Standard distributed algorithms achieve DP typically by perturbing all local subgradients with noise, leading to significantly degenerated utility.
To tackle this issue, we develop a class of private distributed dual averaging (DDA) algorithms, which activates a fraction of nodes to perform optimization.
Such subsampling procedure provably amplifies the DP guarantee, thereby achieving an equivalent level of DP with reduced noise.
We prove that the proposed algorithms have utility loss comparable to centralized private algorithms for both general and strongly convex problems. When removing the noise, our algorithm attains the optimal $\mathcal{O}(1/t)$ convergence for non-smooth stochastic optimization.
Finally, experimental results on two benchmark datasets are given to verify the effectiveness of the proposed algorithms.
\end{abstract}

\end{frontmatter}

\section{Introduction}
\label{sec:intro}

Consider a group of $n$ nodes, where each node $i$ has a local dataset ${D}_i = \{ \xi_{i}^{(1)}, \dots, \xi_{i}^{(q)} \}$ that contains a finite number $q$ of data samples. The nodes are connected via a communication network. 
They aim to collaboratively solve the empirical risk minimization (ERM) problem, where the machine learning models are trained by minimizing the average of empirical prediction loss over known data samples.
Formally, the optimization problem is given by
\begin{equation}\label{OPT}
	\min_{x\in\mathbb{R}^m}  \left\{F(x) := \frac{1}{n}\sum_{i=1}^{n}
 f_i(x)
+h(x) \right\} ,
\end{equation}
where 
$
 f_i(x) = \frac{1}{q} \sum_{j=1}^{q}{l}_i(x,\xi_{i}^{(j)})
$
represents the empirical risk on node $i$,
$l_i(x,\xi)$ is the loss of the model $x$ over the data instance $\xi$, and $h(x)$ is the regularization term shared across the nodes.
{This setup has been commonly considered in machine learning \cite{lian2017can}, where $h(x)$ is used to promote sparsity or model the constraints.}

{As the loss and its gradient in ERM are characterized by data samples, potential privacy issues arise when the datasets are sensitive \cite{bassily2014private}.
In particular, when $l_i$ is the hinge loss, the solution to Problem \eqref{OPT}, i.e., support vector machine (SVM), in its dual form typically discloses data points \cite{bassily2014private}.} Advanced attacks such as input data reconstruction \cite{zhu2020deep} and attribute inference \cite{melis2019exploiting} can extract private information from the gradients. 
To defend privacy attacks, differential privacy (DP) has become prevalent in cryptography and machine learning \cite{dwork2006differential,abadi2016deep},  due to its precise notion and computational simplicity. Informally, DP requires the outcome of an algorithm to remain stable under any possible changes to an individual in the database, and therefore protects individuals from attacks that try to steal the information particular to them.
The DP constraint induces a tradeoff between privacy and utility in learning algorithms \cite{bassily2014private,kifer1private,chaudhuri2012near,wang2017differentially}.

In this work, we are interested in solving Problem \eqref{OPT} while providing rigorous DP guarantee for each data sample in $D:=\cup_{i=1}^nD_i$.

{

}

\subsection{Related work}



For problems without regularization, i.e., $h\equiv 0$, the authors in \cite{huang2015differentially} developed a differentially private distributed gradient descent (DGD) algorithm by perturbing the local output with Laplace noise. Notably, the learning rate is designed to be linearly decaying such that the sensitivity of the algorithm also decreases linearly.
{Then, one can decompose the prescribed DP parameter $\varepsilon$ into a sequence $\{\varepsilon_t\}_{t\geq 1}$, such that $\sum_{\tau=1}^\infty\varepsilon_\tau =\varepsilon$ and the operation at each time instant $t$ can be made $\varepsilon_t$-DP.}
However, such choice of learning rate slows down the convergence dramatically and results in a utility loss in the order of $\mathcal{O}(m/\varepsilon^2)$, where $m$ denotes the dimension of the decision variable.
Under the more reasonable learning rate $\Theta(1/\sqrt{t})$, the utility loss can be improved to $\mathcal{O}\left(\sqrt[4]{{mn^2}/{\varepsilon}}\right)$ \cite{han2016differentially}, where $n$ denotes the number of nodes.
Along this line of research, the authors in \cite{zhu2018differentially,xiong2020privacy,han2021differentially} extended the algorithm to time-varying objective functions, and the authors in \cite{ding2021differentially} advanced the convergence rate to linear based on an additional gradient-tracking scheme.
The authors in \cite{wang2022differentially} developed a distributed algorithm with DP for stochastic aggregative games. The differentially private distributed optimization problem with coupled equality constraints has been studied in \cite{chen2021distributed}.
In these works, however, $\varepsilon$-DP is proved only for each iteration, leading to a cumulative privacy loss of $t\varepsilon$ after $t$ iterations.
{To attenuate the noise effect while ensuring DP, the authors in \cite{vlaski2021graph}  constructed topology-aware noise, with which each node 
perturbs the messages to its neighbors (including itself) with different perturbations whose weighted sum is $0$.}


{For federated learning (FL) with heterogeneous data, the authors in \cite{hu2020personalized} developed a personalized linear model training algorithm with DP. In \cite{noble2022differentially}, general models were considered. In particular, the subsampling of users and local data has been explicitly considered to amplify the DP guarantee and improve the utility.
}

To tackle regularized learning problems, the alternating direction method of multipliers (ADMM) has been used to design distributed algorithms with DP \cite{zhang2016dynamic,zhang2018admm,zhang2018improving}. However, 
an explicit tradeoff analysis between privacy and utility was missing.
The authors in \cite{xiao2021towards} investigated the privacy guarantee produced not only by random noise injection but also by \emph{mixup} \cite{zhang2018mixup}, i.e., a random convex combination of inputs. Approximate DP and advanced composition \cite{kairouz2015composition} were used to keep track of the cumulative privacy loss. The \mbox{privacy--utility} tradeoff in linearized ADMM and DGD were captured by the bound $\mathcal{O}\left({m}/{(\sqrt{n}\varepsilon)}\right)$.

To summarize, existing private distributed optimization algorithms applied to Problem \eqref{OPT} typically require each node to make a gradient query to the local dataset at each time instant. 
{Since the sizes of local datasets are considerably smaller than that of the original dataset, local gradient queries have larger sensitivity parameters than that in centralized settings. 
Therefore, private distributed optimization paradigms in the literature typically employed a larger magnitude of noise to secure the same level of DP, and suffered from relatively low utility.
{Recently, an asynchronous DGD method with DP was developed in \cite{xu2021dp}, which achieved a lower utility loss. The algorithm assumed that each local mini-batch is a subset of data instances uniformly sampled from the overall dataset without replacement, which appears to be restrictive in distributed settings.}


}

\subsection{Contribution}



We develop a class of differentially private distributed dual averaging (DDA) algorithms for solving Problem \eqref{OPT}. 
At each iteration, a fraction of nodes is activated uniformly at random to perform local stochastic subgradient query and local update with perturbed subgradient. Such subsampling procedure provably amplifies the DP guarantee and therefore helps achieve the same level of DP with weaker noise. To ensure a user-defined level of DP, we provide sufficient conditions on the noise variance in Theorem \ref{thm:privacy}, which admits a smaller bound of variance than existing results.




The properties of the proposed algorithms in terms of convergence and the privacy-utility tradeoff are analyzed. First, a non-asymptotic convergence analysis is conducted for dual averaging with inexact oracles under general choices of hyperparameters, and the results are summarized in Theorem \ref{thm_error_bound}. This piece of result illustrates how the lack of global information and the DP noise in private DDA quantitively affect the convergence, which lays the foundation for subsequent analysis. 
Then, we investigate the convergence rate of the non-private (noiseless) version of DDA for both strongly convex and general convex objective functions under two sets of hyperparameters in Corollaries \ref{SC_rate} and \ref{C_rate}, respectively. We remark that Corollary \ref{SC_rate} advances the best known convergence rate of DDA  for nonsmooth stochastic optimization, i.e., $\mathcal{O}(1/\sqrt{t})$, to $\mathcal{O}(1/t)$.
{The key to obtaining the improved rate is the use of a new class of parameters.}

The privacy--utility tradeoff of the proposed algorithm is examined in Corollaries \ref{Cor:SC_loss} and \ref{Cor:C_loss}.
In particular, when the objective function is non-smooth and strongly convex, the utility loss is characterized by $\mathcal{O}(m\iota^2/(q^2\varepsilon^2))$, where $m$, $\iota$, $q$, $\varepsilon$ denote the variable dimension, node sampling ratio, number of samples per node, and DP parameter, respectively. For comparison, we present in {Table \ref{tab:table1}} a comparison of some of the most relevant works.\footnote{Table \ref{tab:table1} presents the dependence on variable dimension $m$, number of nodes $n$, number of samples $q$ per node, sampling ratio $\iota \propto 1/n$,  and $\varepsilon$ for utility loss. The work in \cite{xu2021dp} considered nonconvex problems, and the results are adapted to convex problems for comparison in Table \ref{tab:table1}.} 

{Finally, we verify the effectiveness of the proposed algorithms via distributed SVM on two open-source datasets. Several comparison results are also presented to support our theoretical findings.}


\subsection{Outline}
The rest of the paper is organized as follows. Section \ref{Sec:Pre} introduces some preliminaries. We present our algorithms and their theoretical properties in Section \ref{Sec:Algorithm}, whose proofs are postponed to Section \ref{Sec:Proof}. Some experimental results are given in Section \ref{Sec:Experiment}. Section \ref{Sec:Concl} concludes the paper.

\begin{table*}[!t]

\caption{A Comparison of Some Related Works \label{tab:table1}}
\centering
\begin{tabular}{|c|c|c|c|c|c|c|}
\hline
\multirow{2}{*}{}   & 
\multirow{2}{*}{Privacy} & \multirow{2}{*}{Noise} & {Perturbed} &
\multicolumn{2}{c|}{Utility Upper Bound}  &  {Non-smooth} \\
\cline{5-6} 
& & & Term & Convex& Strongly Convex &  Regularizer\\
\hline
\cite{noble2022differentially} (FL)   & ($\varepsilon$, $\delta$)-DP & Gaussian &   Gradient &$\mathcal{O}\left(\frac{\sqrt{m}}{nq\varepsilon}\right)$& $\mathcal{O}\left(\frac{m}{n^2q^2\varepsilon^2}  \right)$   & No \\
\hline
\cite{huang2015differentially}   & $\varepsilon$-DP & Laplace & Output & --& $\mathcal{O}(\frac{m}{\varepsilon^2})$  & No \\
\hline
\cite{xiao2021towards} (ADMM)   & $(\varepsilon,\delta)$-DP & Gaussian &  Output & $\mathcal{O}\left(  \frac{{n}^{3/2}\varepsilon}{m}\right)$ & -- &  No \\
\hline
\cite{xiao2021towards} (DGD)  & $(\varepsilon,\delta)$-DP & Gaussian & Gradient & ${\mathcal{O}}(\frac{m}{\sqrt{n}\varepsilon})$& --  & No \\
\hline
\cite{xu2021dp}  & $(\varepsilon,\delta)$-DP  & Gaussian & Gradient & ${\mathcal{O}}(\frac{\sqrt{m}}{q\varepsilon})$& --  & No \\
\hline
\textbf{This work}   & ($\varepsilon$, $\delta$)-DP & Gaussian & Gradient & $\mathcal{O}\left(  \frac{\sqrt{m\iota}}{q\varepsilon} \right)$& $\mathcal{O}\left( \frac{m\iota^2}{ q^2\varepsilon^2} \right)$   & Yes \\
\hline
\end{tabular}
\end{table*}

\section{Preliminaries}\label{Sec:Pre}

\subsection{Basic setup}\label{subsec:setup}


We consider the distributed ERM in~\eqref{OPT}, in which $h$ is a closed convex function with non-empty domain $\mbox{dom}(h)$.
Examples of $h(x)$ include  $l_1$-regularization, i.e., $h(x)=\lambda\lVert x \rVert_1, \lambda>0$, and the indicator function of a closed convex set. The regularization term $h$ and the loss functions $l_i$ for all $i=1,\dots,n$ satisfy the following assumptions.

\begin{assumption}
\label{mu-convexity} 
	i) $h(\cdot)$ is a proper closed convex {(strongly convex)} function with modulus $\mu = 0$ {(resp. $\mu>0$)}, i.e., for any $x,y\in dom(h)$,
	\begin{equation*}
		h(\alpha x +(1-\alpha)y) \leq \alpha h(x)+(1-\alpha)h(y)- \frac{\mu \alpha(1-\alpha)}{2}\lVert x-y \rVert^2;
	\end{equation*}
	ii) each $l_i(\cdot, \xi_i)$ is convex on $dom(h)$.
\end{assumption}


When $q=1$, Problem \eqref{OPT} reduces to a deterministic distributed optimization problem. In Problem \eqref{OPT}, the information exchange only occurs between connected nodes.
Similar to existing research \cite{nedic2009distributed,duchi2011dual}, we use a doubly stochastic matrix $W\in[0,1]^{n\times n}$ to encode the network topology and the weights of connected links at time $t$. In particular, its $(i,j)$-th entry, $w_{ij}$, denotes the weight used by $i$ when counting the message from $j$. When $w_{ij}=0$, nodes $i$ and $j$ are disconnected. 

\subsection{Conventional DDA}

The DDA algorithm originally proposed by~\cite{duchi2011dual} can be applied to solve Problem~\eqref{OPT}. In particular, let $d(\cdot)\geq 0$ be a strongly convex function with modulus $1$ on $\mbox{dom}(h)$. 
Each node, starting with $z_i^{(1)}=0$, iteratively generates $\{z_i^{(t)},x_i^{(t)}\}_{t\geq 1}$ according to
\begin{equation}\label{transfored_DA}
	\begin{split}
	x_i^{(t)} = \argmin_{x\in\mathbb{R}^m}\left\{
	\langle z_i^{(t)}, x\rangle + 
t(	h(x) )
	+ \gamma_t d(x)\right\}
\end{split}
\end{equation}
and 
\begin{equation}\label{transfored_DA_z}
			z_i^{(t+1)} = \sum_{j=1}^n w_{ij} \left(z_j^{(t)}+ \hat{g}_j^{(t)}\right),
\end{equation}
where {$\{\gamma_t \}_{t\geq 1}$ is a non-decreasing sequence of parameters}, $w_{ij}$ is the $(i,j)$-th entry of matrix $W$, $\hat{g}_j^{(t)}\in \partial l_j(x_j^{(t)}, \xi_j^{(t)})$ denotes the stochastic subgradient of local loss over $x_j^{(t)}$ with $\xi_j^{(t)}$ uniformly sampled from $D_j$, and $\partial l_j(x_j^{(t)}, \xi_j^{(t)})$ represents the corresponding subdifferential. Throughout the process, each node only passes $z_i$ to its immediate neighbors and updates $x_i$ according to \eqref{transfored_DA}. 
Existing DDA algorithms, when applied to solve Problem \eqref{OPT}, converge as $\mathcal{O}(1/\sqrt{t})$ \cite{duchi2011dual,colin2016gossip}.


{\subsection{Threat model and DP}

In a distributed optimization algorithm, messages bearing information about the local training data are exchanged among the nodes, which leads to privacy risk. 
In this work, we consider the following two types of attackers.

\begin{itemize}
    \item \emph {Honest-but-curious nodes} are assumed to follow the algorithm to perform communication and computation. However, they may record the intermediate results to
infer the sensitive information about the other nodes.

    \item \emph{External eavesdroppers} stealthily listen to the private communications between the nodes.
\end{itemize}

By collecting the confidential messages, the attackers are abele to infer private information about the users \cite{zhu2020deep}. To defend them, we employ tools from DP. Indeed, DP has been recognized as the gold standard in quantifying individual privacy preservation for randomized algorithms. It refers to the property of a randomized algorithm that the presence or absence of an individual in a dataset cannot be distinguished based on the output of the algorithm. Formally, we introduce the following definition of DP for distributed optimization algorithms \cite{zhang2018improving}.




\begin{definition}
\label{def:DPoverZ}
    Consider a communication network, in which each node has its own dataset $D_i$. Let $\{z_i^{(t)}:i=1,\dots,n\}$ denote the set of messages exchanged among the nodes at iteration $t$. A distributed algorithm 
    satisfies $(\varepsilon,\delta)$-DP during $T$ iterations, if for every pair of neighboring datasets $D=\cup_{i=1}^nD_i$ and $D'=\cup_{i=1}^nD'_i$, and for any set of possible outputs $\mathcal{O}$ during $T$ iterations we have
    \begin{equation*}
    \begin{split}
		&Pr[\{ z_i^{(t)},i=1,\dots,n\}_{t=1}^T\in \mathcal{O}|D]\\
		&\leq  e^{\varepsilon} Pr[\{ z_i^{(t)},i=1,\dots,n\}_{t=1}^T\in \mathcal{O}|D']+\delta.
    \end{split}
	\end{equation*}
\end{definition}

}

\section{Differentially private DDA algorithm}\label{Sec:Algorithm}

{In this section, we develop the differentially private DDA algorithm, followed by its privacy-preserving and convergence properties.}

\subsection{Node subsampling in distributed optimization}


{As explained in Section \ref{sec:intro}, parallelized local gradient queries in distributed optimization necessitate stronger noise to achieve DP and therefore deteriorate utility.}
To circumvent this problem, we only activate a random fraction of the nodes at each time instant to perform averaging and local optimization. 
This allows us to amplify the privacy of the algorithm, and thereby achieving the same level of DP with noise weaker than in existing works.

\begin{definition}
\label{participation_rate}
For every $t\geq 1$, an integer number of $n\iota$ nodes are sampled uniformly at random with some $\iota\in(0,1]$. 
\end{definition}


The sampling procedure gives rise to a time-varying stochastic communication network. Slightly adjusted to the notation in Section \ref{subsec:setup}, we let $W^{(t)}\in[0,1]^{n\times n}$ be a {random} gossip matrix at time $t$, where the $(i,j)$-th entry, $w_{ij}^{(t)}$, denotes the weight of the link $(i,j)$ at time $t$. 
Denote by $\mathcal{N}^{(t)}$ and $\mathcal{N}_i^{(t)}:=\{j|j\neq i, w_{ij}^{(t)}>0\}$  the set of activated nodes and the set of $i$'s neighbors at time $t$, respectively. 
{It is worthwhile to point out that $W^{(t)}$ and $\iota$ are dependent. That is, we have $w_{ij}^{(t)} > 0$ for $i\in\mathcal{N}^{(t)}$ and $j\in\mathcal{N}_i^{(t)}$, and $w_{ij}^{(t)}=0$ otherwise.}

For the gossip matrix $W^{(t)}$, we assume the following standard condition \cite{liu2021decentralized}.

\begin{assumption}
\label{spectral-gap}
	For every $t\geq 1$, 	i) $W^{(t)}$ is doubly stochastic\footnote{$W^{(t)}\mathbf{1} = \mathbf{1}$ and $\mathbf{1}^TW^{(t)} = \mathbf{1}^T$ where $\mathbf{1}$ denotes the all-one vector of dimensionality $n$.};
	ii) $W^{(t)}$ is independent of the random events that occur up to time $t-1$; and
	iii) there exists a constant $\beta\in(0,1)$ such that
	\begin{equation}
		\label{eq:sigma-bound}
		\sqrt{\rho\left(\mathbb{E}_{W}\left[{W^{(t)}}^TW^{(t)}\right] - \frac{\mathbf{1}\mathbf{1}^T}{n}\right)} \leq \beta,
	\end{equation}
	where $\rho(\cdot)$ denotes the spectral radius and the expectation $\mathbb{E}[\cdot]$ is taken with respect to the distribution of $W^{(t)}$ at time $t$.
\end{assumption}


\subsection{Private DDA with stochastic subgradient perturbation}

{Next, we introduce a differentially private DDA algorithm presented as Algorithm \ref{DP-DDA}.}

The update for the local dual variable $z_i^{(t)}$ reads
\begin{equation}\label{DP_z}  
	z_i^{(t+1)} =\sum_{j=1}^nw_{ij}^{(t)} \left( z_{j}^{(t)}+ a_t \eta_j^{(t)} \zeta_j^{(t)}\right),
\end{equation}
where $
\zeta_j^{(t)} = \hat{g}_j^{(t)} +\nu_j^{(t)}
$, $\nu_j^{(t)}\sim \mathcal{N}(0, \sigma^2 I)$,
$\hat{g}_j^{(t)}\in \partial l_j(x_j^{(t)}, \xi_j^{(t)})$ with $\xi_j^{(t)}$ uniformly sampled from $D_j$,
$\eta_i^{(t)}$ indicates whether node $i$ performs local update at time $t$, i.e.,
    \begin{equation*}
    \eta_i^{(t)}=\left\{
   \begin{aligned}
   &1,  & \mbox{if} \,\, i \,\, \mbox{active}  \\ 
   &0, &  \mbox{ otherwise}
      \end{aligned}
   \right.
\end{equation*}
and $\{a_t>0\}_{t\geq 1}$ is a sequence of non-decreasing parameters.
The non-decreasing property of $\{a_t\}_{t\geq 1}$ is motivated by that, when the objective exhibits some desirable properties, e.g., strong convexity, assigning heavier weights to fresher subgradients can speed up convergence \cite{lu2018relatively,tao2021gradient}.
In the special case where $a_t=1$, $\eta_i^{(t)}=1$ and $\sigma=0$, \eqref{DP_z} reduces to the conventional update in \eqref{transfored_DA_z}.

Equipped with \eqref{DP_z}, node $i$, active at time $t$, can perform a local computation to derive its estimate about the global optimum: \begin{equation}\label{DP_primal}
	\begin{split}
		{x}_{i}^{(t+1)}
		= \argmin_{x\in\mathbb{R}^m}\left\{ \langle z_{i}^{(t+1)},x \rangle+ \iota A_{t+1}
		h(x)
		+\gamma_{t+1}d(x)\right\} ,
	\end{split}	
\end{equation}
where $\iota$ is defined in Definition \ref{participation_rate}, $A_t= \sum_{\tau=1}^t  a_\tau$ and $\{\gamma_t\}_{t\geq 1}$ is a non-decreasing sequence of positive parameters. 
By convention, we let $A_0=a_0=0$ and $\gamma_0=0$.

\begin{algorithm}[tb]
	\caption{Differentially Private DDA}
	\label{DP-DDA}
	\begin{algorithmic}[1]
		\STATEx {\bfseries Input:} $\mu\geq 0$, $a>0$, a strongly convex function $d$ with modulus $1$ on $\mbox{dom}(h)$, and $T>0$ 
		\STATEx {\bfseries Output:} {$\tilde{x}_i^{(T)}= A_T^{-1} \sum_{t=1}^T a_t x_i^{(t)}$}
		\STATE {\bfseries Initialize:} set $z_i^{(1)} = 0$ and identify $x_i^{(1)}$ according to \eqref{DP_primal} for all $i =1,\dots,n$ 
		\FOR{$t=1,2,\dots, T$}
        \STATEx \emph {\quad For active nodes $i\in\mathcal{N}^{(t)}$}:
		\STATE $\xi_i^{(t)}\sim \text{Uniform}\{1,\dots,q\}$ and $\nu_i{(t)}\sim \mathcal{N}(0, \sigma^2 I)$
        \STATE 
        release $\hat{g}_i^{(t)}+\nu_i^{(t)}$
	\STATE collect $z_j^{(t)}+a_t(\hat{g}_j^{(t)}+\nu_j^{(t)})$ from all nodes $j\in \mathcal{N}_i^{(t)}$
		\STATE update $z_i^{(t+1)}$ by \eqref{DP_z}  
		\STATE compute $x_i^{(t+1)}$ by \eqref{DP_primal}
		\STATEx \emph {\quad For inactive nodes $i\notin\mathcal{N}^{(t)}$}:
		\STATE set $z_i^{(t+1)} = z_i^{(t)}$ and $x_i^{(t+1)} = x_i^{(t)}$
		\ENDFOR
	\end{algorithmic}
\end{algorithm}

For general regularization $h(x)$, the update in \eqref{DP_primal} requires the knowledge of $\iota$. This requirement is necessary due to technical reasons. More precisely, due to node sampling, the term $\langle z_i^{(t+1)}, x \rangle$ in \eqref{DP_primal} serves as a linear approximation of $\iota f_i(x)/n$ rather than $f_i(x)/n$ in standard DDA \cite{duchi2011dual}.
Thus, one scales up $h(x)$ also with $\iota$ in \eqref{DP_primal} in order to solve the original problem in \eqref{OPT}. In the special case where $h(x)$ is the indicator function of a convex set, the knowledge of $\iota$ is not needed since $\iota h(x)\equiv h(x)$.

The overall procedure is summarized in Algorithm \ref{DP-DDA}. 
{Each node $i=1,\dots, n$ initializes $z_i^{(1)}=0$ in Step 1. At each time instant $t$, only active nodes $i \in \mathcal{N}^{(t)}$ update $z_i^{(t+1)}$ and $x_i^{(t+1)}$ by following Steps 3--7. In particular, each active node computes and then perturbs the local stochastic subgradient in Step 3 and 4, respectively, followed by the information exchange with neighboring nodes in Step 5. Then, $z_i^{(t+1)}$ and $x_i^{(t+1)}$ are updated in Steps 6 and 7. For inactive nodes at each time instant $t$, they simply set $z_i^{(t+1)}= z_i^{(t)}$ and $x_i^{(t+1)}= x_i^{(t)}$. }

{\begin{remark}
    There are two common approaches to achieve DP for optimization methods. The first type disturbs the output of a non-private algorithm \cite{zhang2017efficient}, and the second type perturbs the subgradient \cite{bassily2014private,wang2017differentially}.
The former involves recursively estimating the (time-varying) sensitivity of updates \cite{wang2023tailoring}. This makes the propagation of DP noise and its effect on convergence difficult to quantify \cite{wang2023tailoring}.
In this work, we adopt the latter approach in Algorithm \ref{DP-DDA}, where we introduce Gaussian noise to perturb the stochastic subgradient $\hat{g}_{i}$. By leveraging the time-invariant sensitivity of the gradient query, we can effectively conduct both privacy and utility analyses in the presence of non-smooth regularization. It is worth noting that, in this scenario, the step-size scheduling rule allows for control over utility.


\end{remark}

}

\subsection{Privacy analysis}

To establish the privacy-preserving property of Algorithm~\ref{DP-DDA}, we make the following assumption.



\begin{assumption}
\label{L-Lipschitz}
	Each $l_i(\cdot,\xi_i)$ is $L$-Lipschitz, i.e., for any $ x,y\in dom(h)$,
	\begin{equation*}
		\lvert l_i(x,\xi_i) - l_i(y,\xi_i) \rvert \leq L \lVert x-y \rVert.
	\end{equation*} 
\end{assumption}

By Assumption \ref{L-Lipschitz}, we readily have that each $f_i(\cdot)$ is $L$-Lipschitz.
Next, we state the privacy guarantee for Algorithm \ref{DP-DDA}. The proof can be found in Section \ref{thm:DP}.



\begin{thm}
\label{thm:privacy}
Suppose Assumption \ref{L-Lipschitz} is satisfied, and a random fraction of the nodes with ratio $\iota \in(0,1]$ is active at each time instant. Given parameters ${q}$, $\varepsilon\in(0,1]$, and $\delta_0\in(0,1]$, if 
\begin{equation*}
 \sigma^2 \geq  \frac{32\iota^2L^2 T \log(2/\delta_0)}{q^2\varepsilon^2} 
\end{equation*}
and $T \geq 
5q^2\varepsilon^2/(4 \iota^2)$,
then Algorithm \ref{DP-DDA} is $(\varepsilon,\delta)$-DP for some constant $\delta\in(0,1]$.
\end{thm}


\begin{remark}
A few remarks on the results in Theorem \ref{thm:privacy} are in order:
\begin{itemize}
    
    \item[i)] It can be verified from the proof of Theorem \ref{thm:privacy} that Algorithm \ref{DP-DDA} is ($\hat{\varepsilon},\delta$)-DP after $t\leq T$ iterations with 
\begin{equation}\label{privacy_loss_intermediate}
    \hat{\varepsilon} \leq \sqrt{\frac{3t}{5T}}{\varepsilon}+ \frac{t}{5T}{\varepsilon^2}.
\end{equation}

\item[ii)]  Theorem \ref{thm:privacy} emphasizes that, to achieve a prescribed privacy budget during $T$ iterations, the noise variance $\sigma^2$ is related to the DP parameters $(\varepsilon,\delta)$, the Lipschitz constant $L$ of the loss, the number of samples per local dataset, and the iteration number $T$. Notably, the lower bound for variance is weighted by $\iota^2\leq 1$, meaning that the same level of DP can be achieved with reduced noise.
\end{itemize}
\end{remark}


\subsection{Privacy--utility tradeoff}

Next, we perform a non-asymptotic analysis of Algorithm \ref{DP-DDA}, followed by an explicit privacy--utility tradeoff.

Motivated by \cite{duchi2011dual}, we define an auxiliary sequence of variables:
\begin{equation}
	\label{eq:def-y-DP}
	y^{(t+1)}
	= \argmin_{x\in\mathbb{R}^m}\left\{  \langle \overline{z}^{(t+1)},x \rangle+ \iota A_{t+1}h(x)+\gamma_{t+1} d(x)\right\},
\end{equation}
where $\overline{z}^{(t)} = \frac{1}{n}\sum_{i=1}^nz_i^{(t)}$ and $\{z_i^{(t)}: i=1,\dots,n\}_{t\geq 1}$ are generated by Algorithm~\ref{DP-DDA}. The convergence property of $y^{(t)}$ is summarized in Theorem \ref{thm_error_bound}, whose proof is provided in Section \ref{thm:Convergence}.

\begin{thm}
\label{thm_error_bound}
	Suppose Assumptions \ref{mu-convexity}, \ref{spectral-gap}, and \ref{L-Lipschitz} are satisfied. For all $t\geq 1$, we have
	\begin{equation}\label{error_bound}
	\begin{split}
	&A_t \mathbb{E}[F(\tilde{y}^{(t)})-F(x^*)] \\
	&\leq  \frac{\gamma_t}{\iota} d(x^*)  + \sum_{\tau=1}^{t}\frac{ a_\tau^2}{\mu \iota A_\tau+\gamma_\tau} \Big(M+\frac{m  \iota \sigma^2}{2} + \frac{2\sqrt{m\iota }L\sigma}{1-\beta} \Big),
\end{split}
\end{equation}
where $\sigma$ is defined in Theorem \ref{thm:privacy}, $\tilde{y}^{(t)} = A_t^{-1}\sum_{\tau=1}^{t}a_{\tau}y^{(\tau)}$, $M ={\iota L^2}/{2}+ {2\sqrt{\iota}L^2}/{(1-\beta)} $, and the expectation is over the randomness of the algorithm.
\end{thm}


From the error bound in \eqref{error_bound}, we observe that the last two terms are contributed by the noise. How the noise affects the error bound is determined in part by the hyperparameters of the algorithm.
Next, we first investigate the choices of $a_t$ that lead to optimal convergence rates for Algorithm \ref{DP-DDA} with $\sigma=0$; the results for strongly convex and general convex functions are presented in Corollaries \ref{SC_rate} and \ref{C_rate}, whose proofs are given in Appendices \ref{Cor:SC_rate} and \ref{Cor:C_rate}, respectively.

\begin{corollary}
\label{SC_rate}
	{Suppose Assumptions \ref{spectral-gap} and \ref{L-Lipschitz} are satisfied. In addition, Assumption \ref{mu-convexity} holds with $\mu>0$, i.e., $h(x)$ is $\mu$-strongly convex.
 If $\sigma=0$,} 
	\begin{equation*}
		a_t=t\ \ \mbox{and} \ \  \gamma_t= 0,
	\end{equation*}
 then for all $t\geq 1$, and $i=1,\dots,n$, we have
	\begin{equation}\label{rate_mu_convex}
		\begin{split}
		\frac{1}{n}	\sum_{i=1}^n	\mathbb{E}[\lVert  \tilde{x}_i^{(t)}-x^*\rVert^2] \leq \frac{16}{t+1}\left(  \frac{L^2 (\log t+1)}{\mu^2\iota(1-\beta)^2t}+ \frac{M }{\mu^2 \iota }\right),
		\end{split}
	\end{equation}
where $\tilde{x}_i^{(t)}= A_t^{-1}\sum_{\tau=1}^t a_\tau x_i^{(\tau)}$ and $M$ is a positive constant given in Theorem \ref{thm_error_bound}.
\end{corollary}

{\begin{remark}
Corollary \ref{SC_rate} indicates that the non-private version of Algorithm \ref{DP-DDA}, i.e., $\sigma=0$, attains the optimal convergence rate $\mathcal{O}(1/t)$ when Problem \eqref{OPT} is strongly convex. Compared to the algorithm in \cite{yuan2018optimal}, where the authors focused on constrained problems, the proposed algorithm handles general non-smooth regularizers. Furthermore, the results can be extended to the case where each $f_i(x)$ but not $h(x)$ is strongly convex by following a similar idea in \cite{liu2021decentralized}.

\end{remark}
}

\begin{corollary}
\label{C_rate}
		{Suppose Assumptions \ref{spectral-gap} and \ref{L-Lipschitz} are satisfied. In addition, Assumption \ref{mu-convexity} holds with $\mu=0$, i.e., $h(x)$ is general convex. If $\sigma=0$, }
	\begin{equation*}
		a_t=1 \ \ \mbox{and} \ \  \gamma_t= \gamma \sqrt{t},
	\end{equation*}
	then for all $t\geq 1$, and $i=1,\dots,n$, we have
	\begin{equation}\label{objective_error_y}
		\mathbb{E}\left[ F(\tilde{y}^{(t)})-F(x^*) \right]\leq\frac{d(x^*) +2\iota M}{\iota\gamma\sqrt{t}},
	\end{equation}
where $\tilde{y}^{(t)}=t^{-1}\sum_{\tau=1}^{t}y^{(\tau)}$ and $M$ is a positive constant given in Theorem \ref{thm_error_bound}.
 In addition, for all $t\geq 1$, and $i=1,\dots,n$, we have
	 \begin{equation}
	 \label{consensus_error_convex}
    \frac{1}{n}\sum_{i=1}^n\mathbb{E}[\lVert \tilde{x}_i^{(t)} - \tilde{y}^{(t)}\rVert]  \leq  \frac{2L\sqrt{\iota}}{\gamma(1-\beta)\sqrt{t}},
\end{equation}
	where $\tilde{x}_{i}^{(t)}=t^{-1}\sum_{\tau=1}^{t}{x}_{i}^{(\tau)}$.
\end{corollary}

Under the same hyperparameters, we study the privacy--utlity tradeoff of Algorithm \ref{DP-DDA} with $\sigma\neq 0$ for strongly convex and general convex functions in Corollaries \ref{Cor:SC_loss} and \ref{Cor:C_loss}, whose proofs are presented in Appendices \ref{App:SC_loss} and \ref{App:C_loss}, respectively.
\begin{corollary}
\label{Cor:SC_loss}
{Suppose Assumptions \ref{spectral-gap} and \ref{L-Lipschitz} are satisfied. In addition, Assumption \ref{mu-convexity} holds with $\mu>0$, i.e., $h(x)$ is $\mu$-strongly convex. If}	
	\begin{equation*}
		a_t={t},\ \  \ \gamma_t= 0,
	\end{equation*} 
	and $\beta \in (1-\sqrt{1/e}, 1]$,
then for $T = \mathcal{O}\left(\frac{q^2\varepsilon^2}{\iota^3  (1-\beta)^2m\log(1/\delta_0)}\right)$ and $i=1,\dots, n$:
\begin{equation}\label{ULoss:SC}
\begin{split}
   \frac{1}{n}\sum_{i=1}^n \mathbb{E}\left[\lVert  \tilde{x}_i^{(T)} -  x^* \rVert^2 \right] \leq {\mathcal{O}}\left( \frac{m\iota^2 L^2 \log(1/\delta_0) }{\mu^2 q^2\varepsilon^2}  \right),
    \end{split}
\end{equation}
where $\tilde{x}_i^{(T)}= A_T^{-1} \sum_{t=1}^T a_t x_i^{(t)}$, and $\delta_0$ is defined in Theorem \ref{thm:privacy}.
\end{corollary}

\begin{corollary}
\label{Cor:C_loss}
{Suppose Assumptions \ref{spectral-gap} and \ref{L-Lipschitz} are satisfied. In addition, Assumption \ref{mu-convexity} holds with $\mu=0$, i.e., $h(x)$ is general convex. If}
\begin{equation}\label{parameter_choices_C_DP}
	a_t=1 , \  \  \gamma_t= \gamma \sqrt{t},   
\end{equation}
and $ \iota \leq 1-\beta $,
then the following holds for $T=\mathcal{O}\left(\frac{q^2\varepsilon^2}{\iota^3(1-\beta)^2{m\log(1/\delta_0)}}\right)$:
\begin{equation}\label{ULoss:C}
\begin{split}
    \mathbb{E}[F(\tilde{y}^{(T)})-F(x^*)] \leq \mathcal{O}\left(  \frac{\left(L^2 +d(x^*)\right)\sqrt{m\iota \log(1/\delta_0)}}{\gamma q\varepsilon } \right)
\end{split}
\end{equation}
and 
\begin{equation*}
\begin{split}
		\frac{1}{n}\sum_{i=1}^n\mathbb{E}[\lVert \tilde{x}_i^{(T)} - \tilde{y}^{(T)}\rVert]   	\leq  \frac{L\sqrt{m\iota \log(1/\delta_0)}}{\gamma q\varepsilon}
\end{split}
\end{equation*}
for $i=1,\dots,n$,
where $\tilde{y}^{(T)}=\frac{1}{T}\sum_{t=1}^{T}{y}^{(t)}$, $\tilde{x}_i^{(T)}=\frac{1}{T}\sum_{t=1}^{T}{x}_i^{(t)}$, and $\delta_0$ is defined in Theorem \ref{thm:privacy}.
\end{corollary}

Corollaries \ref{Cor:SC_loss} and \ref{Cor:C_loss} highlight that the sampling procedure lowers down the utility loss for both strongly convex and general convex problems. In particular, the utility loss in the strongly convex case becomes $\iota^2\approx 1/n^2$ times smaller than that without sampling. For general convex problems, the utility loss is  $\sqrt{\iota}$ times smaller.
They also suggest that the number of iterations increases in order to achieve a lower utility loss.

\section{Proofs of main results}\label{Sec:Proof}

This section presents the proof of Theorems \ref{thm:privacy} and \ref{thm_error_bound}.

\subsection{Proof of Theorem \ref{thm:privacy}\label{thm:DP}}

We start by introducing some useful properties of DP \cite{dwork2006differential,kairouz2015composition,girgis2021shuffled}.

\begin{lemma}[Gaussian Mechanism]\label{Gaussian_to_zCDP}
	Consider the Gaussian mechanism for answering the query $r:\mathcal{D}\rightarrow \mathbb{R}^m$:
		\begin{equation*}
		\mathcal{\mathcal{M}} = r({D}) +\nu,
	\end{equation*} 
where $D\in\mathcal{D}$, $\nu\sim \mathcal{N}(0, \sigma^2 I)$.
The mechanism $\mathcal{M}$ is ($\sqrt{2\log(2/\delta)}\Delta/\sigma,\delta$)-DP where $\Delta$ denotes the sensitivity of $r$, i.e.,
	$
	\Delta  = \sup_{{D}, {D}'}\lVert  r( {D})-r( {D}') \rVert.
	$
\end{lemma}

{Recall from Algorithm \ref{DP-DDA} that, at each iteration, $n\iota$ nodes are sampled from $n$ nodes at random, and each activated node randomly selects a data sample from $q$ instances to compute stochastic gradients. Although such a sub-sampling is not uniform, i.e., the subsets of $n\iota$ data samples are not necessarily chosen with equal probability, it still helps amplify the privacy \cite[Lemma 10]{girgis2021shuffled}.}


\begin{lemma}[Privacy Amplification by Subsampling]\label{Privacy_amplification}
    Suppose $\mathcal{M}$ is an $(\varepsilon,\delta)$-DP mechanism. 
    Let $\text{samp}_{r_1,r_2}: \mathcal{U}^{r_1}\rightarrow \mathcal{U}^{r_2}$ be the subsampling operation that takes a dataset belonging to $ \mathcal{U}^{r_1}$ as input and selects uniformly at random a subset of $r_2\leq r_1$ elements from the input dataset. Then, the mechanism
\begin{equation*}
    \mathcal{A}(\mathcal{D}) = \text{samp}_{n,n\iota} ( \mathcal{G}_1,\dots,\mathcal{G}_n  )
\end{equation*}
where $\mathcal{G}_i=\text{samp}_{q,1}(\mathcal{M}(\hat{g}_i(\xi_i^{(1)})),\dots, \mathcal{M}(\hat{g}_i(\xi_i^{(q)})))$ and $\hat{g}_i(\xi)$ is a subgradient of $l_i$ evaluated over data point $\xi$
is $(\ln(1+\iota(e^\varepsilon-1)/q),\iota\delta/q)$-DP,
\end{lemma}


\begin{lemma}[Composition of DP]\label{composition}
	Given $T$ randomized algorithms $\mathcal{A}_1,\dots,\mathcal{A}_\tau, \dots,\mathcal{A}_T: \mathcal{D}\rightarrow\mathcal{R}$, each of which is ($\varepsilon_i,\delta_i$)-DP with $\varepsilon_i\in(0,0.9]$ and $\delta_i\in(0,1]$. Then $\mathcal{A}:\mathcal{D}\rightarrow \mathcal{R}^t$ with $\mathcal{A}(\mathcal{D})=(\mathcal{A}_1(\mathcal{D}),\dots,\mathcal{A}_t(\mathcal{D}))$ is $(\tilde{\varepsilon},\tilde{\delta})$-DP with
	\begin{equation*}
	    \tilde{\varepsilon}= \sqrt{\sum_{i=1}^k 2\varepsilon_i^2 \log(e+\frac{\sqrt{\sum_{i=1}^k\varepsilon_i^2}}{\delta'})} + \sum_{i=1}^k \varepsilon_i^2
	\end{equation*}
	and 
	\begin{equation*}
	    \tilde{\delta} = 1-(1-\delta')\prod_{i=1}^k(1-\delta_i)
	\end{equation*}
	for any $\delta'\in(0,1]$.
\end{lemma}

{
\begin{lemma}[Post-Processing]\label{lem:post-processing}
    Given a randomized algorithm $\mathcal{A}$ that is $(\varepsilon,\delta)$-DP. For arbitrary mapping $p$ from the set of possible outputs of $\mathcal{A}$ to an arbitrary set, $p(\mathcal{A}(\cdot))$ is $(\varepsilon,\delta)$-DP.
\end{lemma}}


We are now in a position to prove Theorem \ref{thm:privacy}.

 \textbf{DP at each time $t$}:
We begin by noting that the subgradient perturbation procedure at time $t$, denoted by $\mathcal{M}_t$, is a Gaussian mechanism whose sensitivity, by Assumption \ref{L-Lipschitz}, is
$
	\Delta\leq { 2L}.
$
Based on Lemma \ref{Gaussian_to_zCDP}, $\mathcal{M}_t$ is ($\varepsilon_t, \delta_0$)-DP with
\begin{equation*}
    \varepsilon_t = \frac{2L\sqrt{2\log(2/\delta_0)}}{   \sigma}
\end{equation*}
for any $\delta_0\in [0,1]$. {Due to the conditions on $\sigma$ and $T$}, we obtain
\begin{equation}\label{bound_epsi_t}
\begin{split}
    \varepsilon_t^2 =   \frac{8 L^2 \log (2/\delta_0) }{ \sigma^2} = \frac{ q^2 \varepsilon^2 }{4\iota^2T } \leq 0.2.
        \end{split}
\end{equation}
Denote by $\mathcal{A}_t$ the composition of $\mathcal{M}_t$ and the subsampling procedure.
Upon using Lemma \ref{Privacy_amplification} and \eqref{bound_epsi_t}, we obtain that $\mathcal{A}_t$ is $(\varepsilon_t', \iota \delta_0/q)$-DP with
\begin{equation*}
         \varepsilon_t' =  \frac{2\iota\varepsilon_t}{q} \geq \frac{\iota(e^{\varepsilon_t}-1)}{q} \geq \ln\left(1+ \frac{\iota(e^{\varepsilon_t}-1)}{q}\right).
\end{equation*}
In addition, because of \eqref{bound_epsi_t}, we get
$
    \varepsilon_t' = 2\iota \varepsilon_t/q \leq 2\varepsilon_t \leq 0.9
$
and
\begin{equation}\label{bound_sum_epsi_p}
\begin{split}
    &\sum_{t=1}^T{\varepsilon_t'}^2 =  \frac{4\iota^2}{q^2}\sum_{t=1}^T \frac{8 L^2 \log (2/\delta_0) }{\sigma^2} = \sum_{t=1}^T\frac{ \varepsilon^2 }{ T} \leq 1.
        \end{split}
\end{equation}



 \textbf{DP after $T$ iterations}:
Consider the composition of $\mathcal{A}_1,\dots,\mathcal{A}_\tau, \dots,\mathcal{A}_T$, denoted by $\mathcal{A}$. Based on the advanced composition rule for DP in Lemma \ref{composition}, we obtain $\mathcal{A}$ is ($\tilde{\varepsilon},\tilde{\delta}$)-DP with
\begin{equation*}
\begin{split}
\tilde{\varepsilon} =\sqrt{2\sum_{t=1}^T {\varepsilon_t'}^2 \log(e+\frac{\sqrt{\sum_{t=1}^T{\varepsilon_t'}^2}}{\delta'})}  + \sum_{t=1}^T {\varepsilon_t'}^2
\end{split}
\end{equation*}
and
$
    \tilde{\delta} = 1 - (1-\delta')(1-\iota\delta_0)^T 
$
for any $\delta'\in(0,1]$.
Furthermore, there holds
\begin{equation*}
\begin{split}
    &\tilde{\varepsilon}\leq \sqrt{2\sum_{t=1}^T {\varepsilon_t'}^2 \log(e+\frac{\sqrt{\sum_{t=1}^T{\varepsilon_t'}^2}}{\delta'})} + \frac{1}{5}\varepsilon^2 \\
    & \leq  \sqrt{3\sum_{t=1}^T{\varepsilon_t'}^2} + \frac{1}{5}\varepsilon = \sqrt{\frac{3}{5}\varepsilon^2} + \frac{1}{5}\varepsilon \\
    &\leq \varepsilon,
    \end{split}
\end{equation*}
where we fix $\delta'=\sqrt{\sum_{t=1}^T{\varepsilon_t'}^2}\leq 1$ and use \eqref{bound_sum_epsi_p} to get the second inequality. By setting $\delta= \tilde{\delta}$, we have that $\mathcal{A}$ is $(\varepsilon,\delta)$-DP.

{
\textbf{DP after postprocessing}:
The intermediate results $\{ {\bf z}^{(\tau)} \}_{\tau=1}^T$ are computed based on the output of $\mathcal{A}$, i.e., perturbed subgradients.
By the post-processing property of DP in Lemma \ref{lem:post-processing}, Algorithm \ref{DP-DDA} also satisfies ($\varepsilon,\delta$)-DP specified in Definition \ref{def:DPoverZ}.}

\subsection{Proof of Theorem \ref{thm_error_bound}\label{thm:Convergence}}

Before proving Theorem \ref{thm_error_bound}, we present two useful lemmas whose proofs are given in Appendices \ref{Appen:Lem4} and \ref{Appen:Lem5}.

\begin{lemma}\label{consensus_error}
 	 For the sequence $\{x_i^{(t)}: i =1,\dots, n\}_{t\geq 1}$ generated by Algorithm~\ref{DP-DDA} and the auxiliary sequence $\{y^{(t)}\}_{t\geq 1}$ defined in~\eqref{eq:def-y-DP}, one has that for all $t\geq 1$ and $i=1,\dots,n$,
	 \begin{equation}\label{error_consensus}
	 \frac{1}{n}\sum_{i=1}^n\mathbb{E}\left[	\lVert   x_i^{(t)}-y^{(t)}\rVert \right] \leq \frac{a_t(L+\sqrt{m}\sigma)\sqrt{\iota}}{(1-\beta)(\mu \iota A_t+\gamma_t)}
	 \end{equation}
 and 
 	\begin{equation}\label{MSE_consensus_expectation}
 	 \frac{1}{n}\sum_{i=1}^n\mathbb{E}\left[	\lVert   x_i^{(t)}-y^{(t)}\rVert^2 \right] \leq \frac{a_t^2(L^2+m\sigma^2){\iota}}{(1-\beta)^2(\mu \iota A_t+\gamma_t)^2}.
 \end{equation}
\end{lemma}

\begin{lemma} \label{auxiliary_sequence}
	For all $t\geq 1$, we have
	\begin{equation}\label{chap2:dual_averaging_inequality}
		\begin{split}
		& \frac{1}{n}\sum_{i=1}^n \sum_{\tau=1}^{t} a_{\tau}  \left( \frac{1}{\iota } \left\langle  \eta_i^{(\tau)}\zeta_i^{(\tau)}, {y}^{(\tau)}-x^* \right\rangle   + h(y^{(\tau)}) -h(x^*)    \right) \\
		&\leq \frac{1}{2} \sum_{\tau=1}^{t}\frac{a_\tau^2}{\iota \left(\mu \iota A_\tau+\gamma_{\tau} \right)} \left\lVert \frac{1}{n}\sum_{i=1}^n\eta_i^{(\tau)} \zeta_i^{(\tau)}\right\rVert^2+\frac{\gamma_t}{\iota} d(x^*).
	\end{split}
	\end{equation} 
\end{lemma}

Now we are ready to prove Theorem \ref{thm_error_bound}.
Upon using $A_t = \sum_{\tau=1}^{t}a_\tau$, convexity of $f:={n}^{-1}\sum_{i=1}^{n}f_i$, and $L$-Lipschitz continuity of each $f_j$, we have
\begin{equation}\label{ergodic_and_non}
	\begin{split}
		&A_t\left(f(\tilde{y}^{(t)})-f(x^*)\right)\\
		& \leq   \sum_{\tau=1}^{t}  a_{\tau}\left(f(y^{(\tau)})-f(x^*)\right) \\
		& = \frac{1}{n}\sum_{\tau =1}^{t} \sum_{j=1}^{n}  a_{\tau}\big(  f_j(y^{(\tau)}) -f_j(x_j^{(\tau)})  +f_j(x_j^{(\tau)}) -f_j(x^*)  \big) \\
		& \leq  \frac{1}{n}\sum_{\tau =1}^t\sum_{j=1}^{n}    a_{\tau} L\lVert  y^{(\tau)} - x_j^{(\tau)}\rVert +\frac{1}{n} \sum_{\tau=1}^{t}{\phi^{(\tau)}} , 
	\end{split}
\end{equation}
where $ \phi^{(\tau)}=\sum_{j=1}^{n} a_\tau\left( f_j(x_j^{(\tau)}) -f_j(x^*)\right)$.
Further using convexity of $f_j, j=1,\dots, n$,
we have
\begin{equation}\label{convexity_2}
	\begin{split}
		&	\phi^{(\tau)}\leq    \sum_{j=1}^{n}  a_\tau\left \langle  g_j^{(\tau)}, x_{j}^{(\tau)}-x^* \right \rangle \\
				& = n  a_\tau \left \langle  \overline{\zeta}^{(\tau)}, y^{(\tau)}-x^* \right \rangle  +	 \sum_{j=1}^{n} a_\tau \Big(\left \langle  g_j^{(\tau)}, x_{j}^{(\tau)}-y^{(\tau)} \right \rangle    \\
				&  \quad + \left \langle g_j^{(\tau)}-\zeta_j^{(\tau)}, y^{(\tau)}-x^* \right \rangle \Big),
	\end{split}
\end{equation}
where 
$\overline{{\zeta}}^{(t)} = n^{-1}\sum_{i=1}^{n}\zeta_i^{(t)}$.
Therefore, 
\begin{equation}\label{ergodic_non_expecations}
	\begin{split}
				&A_t	\left( F(\tilde{y}^{(t)})-F(x^*)\right) \\
		&\leq A_t	\left( f(\tilde{y}^{(t)})-f(x^*)\right) + \sum_{\tau=1}^{t}  a_\tau \left (h(y^{(\tau)}) -h(x^*)\right) \\
		& \leq \sum_{\tau=1}^{t} a_\tau \left(\frac{2}{n}\sum_{j=1}^{n} L \lVert x_j^{(\tau)} - y^{(\tau)} \rVert + h(y^{(\tau)}) - h(x^*) \right)  \\
		& \quad  +   \sum_{\tau=1}^ta_\tau \left \langle   \overline{\zeta}^{(\tau)}, y^{(\tau)}-x^* \right \rangle  \\
		& \quad + \frac{1}{n}\sum_{\tau=1}^{t}\sum_{j=1}^{n} a_{\tau}  \left \langle g_j^{(\tau)}-\zeta_j^{(\tau)},y^{(\tau)}-x^* \right \rangle ,
	\end{split}
\end{equation}
where we use
\begin{equation*}
 \sum_{\tau=1}^{t} a_\tau\left (h(x^*)-h(y^{(\tau)}) \right) \leq A_t\left( h(x^*)- h(\tilde{y}^{(t)}) \right)
\end{equation*}
and $F= f+h$ in the first inequality, and use 
\eqref{ergodic_and_non}, \eqref{convexity_2} in the last inequality. 
Due to uniform node sampling with probability $\iota$, we have
\begin{equation*}
    \frac{1}{n}\mathbb{E}_\tau \left[\sum_{i\in\mathcal{N}^{(t)}}\zeta_i^{(\tau)}\right] = \frac{\iota}{n}  \sum_{i=1}^n\mathbb{E}_\tau \left[ \zeta_i^{(\tau)} \right] =   \iota  \mathbb{E}_\tau \left[ \overline{\zeta}^{(\tau)} \right],
\end{equation*}
where $\mathbb{E}_\tau$ denotes expectation conditioned on $\{ x_i^{(\tau)}, i=1,\dots,n \}$. Therefore, by putting the conditioned expectation on \eqref{ergodic_non_expecations} and using the law of total expectation, we obtain
\begin{equation*}
	\begin{split}
				&A_t	\mathbb{E}\left[ F(\tilde{y}^{(t)})-F(x^*)\right] \\
		& \leq \sum_{\tau=1}^{t} a_\tau \left(\frac{2}{n}\sum_{j=1}^{n} L \mathbb{E}[\lVert x_j^{(\tau)} - y^{(\tau)} \rVert] + \mathbb{E}[h(y^{(\tau)}) - h(x^*)] \right)  \\
		& \quad  +  \frac{1}{\iota} \sum_{\tau=1}^ta_\tau \mathbb{E}\left[\left \langle  \frac{1}{n}\sum_{i=1}^n \eta_i^{(\tau)}\zeta_i^{(\tau)}, y^{(\tau)}-x^* \right \rangle \right] \\
		& \quad + \frac{1}{n}\sum_{\tau=1}^{t}\sum_{j=1}^{n} a_{\tau}  \mathbb{E}\left[\left \langle g_j^{(\tau)}-\zeta_j^{(\tau)},y^{(\tau)}-x^* \right \rangle\right].  
	\end{split}
\end{equation*}
Since $\nu_j^{(\tau)}$ and $\hat{g}_j^{(\tau)}$ are independent of $y^{(\tau)}$ and $\mathbb{E}[\hat{g}_j^{(\tau)}]={g}_j^{(\tau)} $, we have
\begin{equation}\label{zero_error}
	\begin{split}
		&\mathbb{E}\left[  \left\langle g_j^{(\tau)}-\zeta_j^{(\tau)}, y^{(\tau)}-x^* \right\rangle\right]   \\
		& = \mathbb{E}\left[  \left\langle g_j^{(\tau)}-\hat{g}_j^{(\tau)}, y^{(\tau)}-x^* \right\rangle\right] + \mathbb{E}\left[  \left\langle -\nu_j^{(\tau)}, y^{(\tau)}-x^* \right\rangle\right]   \\
		& =  0.
	\end{split}
\end{equation}
Therefore, we obtain from Lemma \ref{auxiliary_sequence} that
\begin{equation}\label{convergence_last}
	\begin{split}
		&A_t\mathbb{E}[F(\tilde{y}^{(t)})-F(x^*)] \\
		&\leq   \frac{1}{2} \sum_{\tau=1}^{t}\frac{a_\tau^2}{\iota \left(\mu \iota A_\tau+\gamma_{\tau}\right)} \mathbb{E}\left[\left\lVert \frac{1}{n} \sum_{i=1}^n \eta_i^{(\tau)} \zeta_i^{(\tau)}\right\rVert^2 \right] \\
		& \quad	+ \frac{\gamma_t}{\iota} d(x^*)
		+{\frac{2}{n}\sum_{\tau=1}^{t} \sum_{j=1}^{n}L a_{\tau}\mathbb{E}\left[\lVert x_j^{(\tau)}- y^{(\tau)}\rVert\right]}.
	\end{split}
\end{equation}
Furthermore, we have
\begin{equation*}
	\begin{split}
	\left\lVert \frac{1}{n} \sum_{i=1}^n \eta_i^{(\tau)} \zeta_i^{(\tau)}\right\rVert^2  = \frac{1}{n^2}  	\left\lVert \sum_{i\in\mathcal{N}^{(t)}}^n  \zeta_i^{(\tau)}\right\rVert^2 \leq \frac{\iota}{n}  \sum_{i\in\mathcal{N}^{(t)}}     \left \lVert  {\zeta}_i^{(\tau)}\right \rVert^2 .
 \end{split}
\end{equation*}
Since
\begin{equation*}
    \mathbb{E}_\tau\left[ \sum_{i\in\mathcal{N}^{(t)}}     \left \lVert  {\zeta}_i^{(\tau)}\right \rVert^2 \right] =  \iota \sum_{i=1}^n  \mathbb{E}_\tau \left[ \left \lVert  {\zeta}_i^{(\tau)}\right \rVert^2 \right],
\end{equation*}
we remove the conditioning based on the law of total expectation to obtain
\begin{equation*}
\begin{split}
    &\mathbb{E}\left[\left\lVert \frac{1}{n} \sum_{i=1}^n \eta_i^{(\tau)} \zeta_i^{(\tau)}\right\rVert^2 \right] \leq      \frac{\iota^2}{n}\sum_{i=1}^n  \mathbb{E} \left[ \left \lVert  {\zeta}_i^{(\tau)}\right \rVert^2 \right]\\
    &\leq \frac{\iota^2}{n}\sum_{i=1}^n  \left( \mathbb{E} \left[ \left \lVert  {\hat{g}}_i^{(\tau)}\right \rVert^2 \right]+\mathbb{E} \left[ \left \lVert  \nu_i^{(\tau)}\right \rVert^2 \right]\right)\\
    &\leq  \iota^2(L^2+m\sigma^2),
\end{split}
\end{equation*}
where we use the fact that $\nu_i^{(\tau)}$ is independent of $\hat{g}_i^{(\tau)}$.
By plugging the above into \eqref{convergence_last} and using Lemma \ref{consensus_error}, we arrive at \eqref{error_bound} as desired.

\section{Experiments}\label{Sec:Experiment}

In this section, we present experimental results of the proposed algorithms.

\subsection{Setup}
We use the benchmark datasets \emph{epsilon} \cite{sonnenburg2006large} and \emph{rcv1} \cite{lewis2004rcv1} in the experiments. Some information about the datasets is given in Table \ref{tab:dataset}. We randomly assign the data samples evenly among the $n=20$ working nodes.
The working nodes aim to solve the following regularized SVM problem:

\begin{equation}\label{eq:svm_l1}
	\min_{x} \left\{F(x) =\frac{1}{n}\sum_{i=1}^{n} f_i(x) + h(x)\right\},
\end{equation}
where 
\begin{equation}\label{eq:svm}
	f_i(x) =\frac{1}{q} \sum_{j=1}^{q} \max \left\{ 0, 1- y_i^{(j)}\left\langle  C_i^{(j)}, x \right\rangle \right\},
\end{equation}
$\{C_i^{(j)}, y_i^{(j)}\}_{j=1}^{q} := D_i$ are data samples private to node $i$. In the experiment, we consider two choices of the regularizer, i.e., $h(x)=\phi\lVert x \rVert_1$ and $h(x)=\mu\lVert x \rVert_2^2$ where $\phi>0$ and $\mu>0$ will be specified later.

\begin{table}[!t]
\caption{Properties of the datasets\label{tab:dataset}}
\centering
\begin{tabular}{|c|c|c|}
\hline
Datasets & \# of samples & \# of features\\
\hline
\emph{epsilon} & 400000 & 2000 \\
\hline
\emph{rcv1} & 677399 & 47236 \\
\hline
\end{tabular}
\end{table}


Throughout the experiments, we consider a complete graph with $n=20$ nodes as the supergraph. Based on it, we consider two edge sampling strategies, that is, $1$ or $2$ edges are sampled uniformly at random from the set of all edges at each time instant.  
The corresponding gossip matrices are created with Metropolis weights \cite{xiao2007distributed}.


Some common parameters used in the two sets of experiments are introduced in the following.
For the parameters of DP, we consider
$ 
    \varepsilon \in \{0.2, 0.4, 0.6, 0.8, 1\}
$ and $\delta_0 = 0.01$.
The random noises in these two cases are generated accordingly based on Theorem \ref{thm:privacy}.
The convergence performance of the algorithm is captured by suboptimality, i.e., $F({n}^{-1}\sum_{i=1}^n\tilde{x}_i^{(t)})-F(x^*)$, versus the number of iterations, where the ground truth is obtained by the optimizer {SGDClassifier} from scikit-learn \cite{scikit-learn}. 






\subsection{Results for $l_2$-regularized SVM}

Set
$
h(x) = \mu\lVert x \rVert^2/2
$
with $\mu = 0.0005$. Since the problem is strongly convex, we set $a_t=t$ and $\gamma_t=20$. 

\begin{figure}[htbp]
\centering
\includegraphics[width=3.3in]{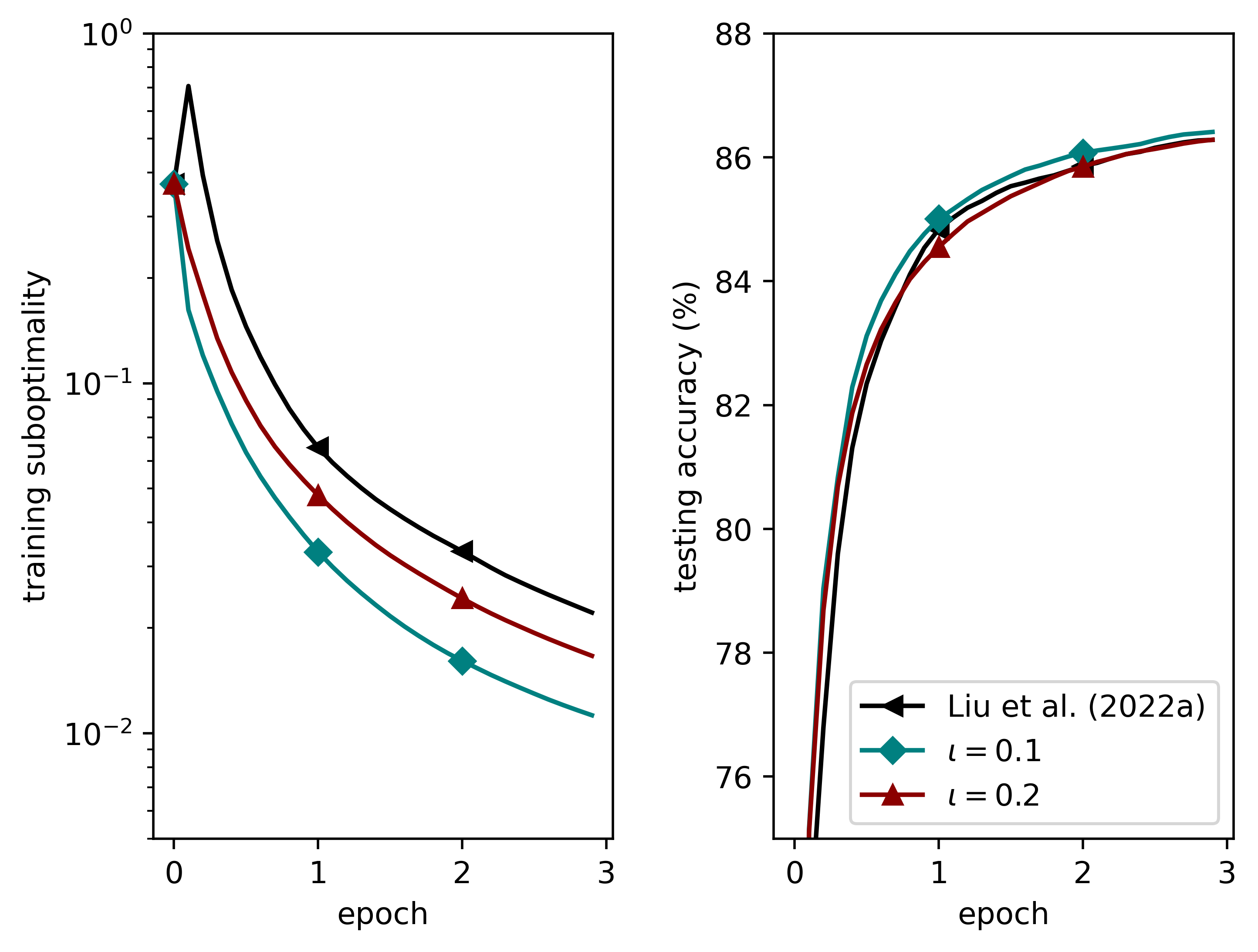}
\caption{Performance comparison between Algorithm \ref{DP-DDA} and \cite{LIU202243} for $l_2$-regularized SVM with $\varepsilon=0.8$.}
\label{l2_subopt}
\end{figure}

We set $\varepsilon=0.8$ and compare the convergence performance between \cite{LIU202243} and Algorithm \ref{DP-DDA} under different choices of $\iota\in\{0.1,0.2 \}$. Fig. \ref{l2_subopt} shows that Algorithm \ref{DP-DDA} with both choices of $\iota$ outperform \cite{LIU202243} in terms of convergence speed and model accuracy. Furthermore, the use of larger $\iota$ in Algorithm \ref{DP-DDA} leads to higher utility loss, which verifies Corollary \ref{Cor:SC_loss}. 
We observe that selecting a higher number of sampled nodes at each step leads to improved network connectivity as well as increased noise. The findings from Fig. \ref{l2_subopt} indicate that, in this specific example, the impact of increased noise on convergence performance may outweigh the benefits of enhanced connectivity.

\begin{figure}[htbp]
\centering
\includegraphics[width=3.3in]{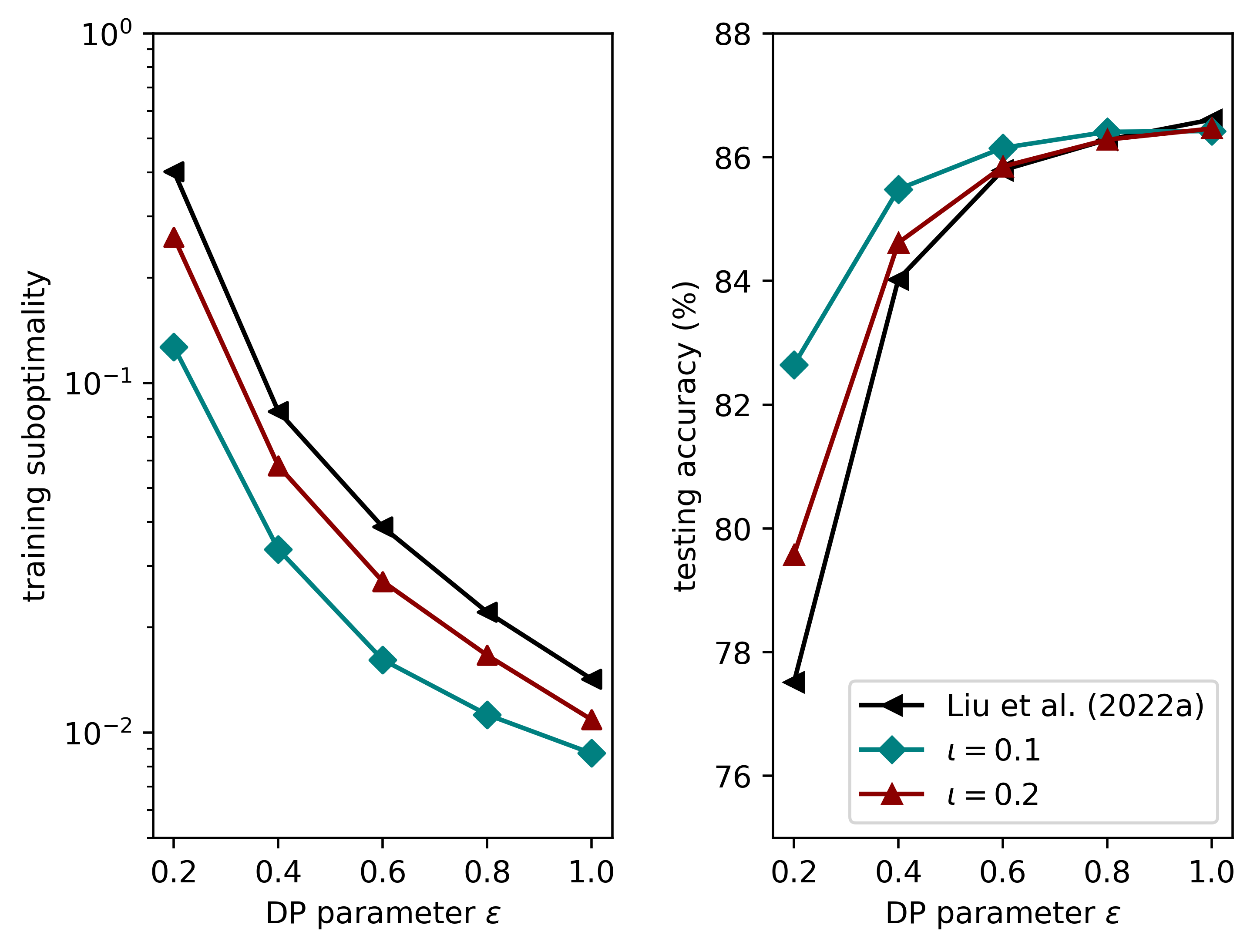}
\caption{Privacy--utility tradeoff in $l_2$-regularized SVM. The suboptimality and accuracy are evaluated after $3$-epoch training.}
\label{l2_nodes_subopt}
\end{figure}

Next, we examine the performance of \cite{LIU202243} and Algorithm \ref{DP-DDA} under a set of DP parameters. The result in Fig. \ref{l2_nodes_subopt} illustrates that 
increasing the value of $\varepsilon$--indicating a less stringent privacy requirement--results in decreased utility loss across all the methods.
This can be attributed to the fact that a smaller value of $\varepsilon$ corresponds to a more stringent differential privacy (DP) constraint, necessitating a stronger noise to perturb the subgradient.
In addition, the performance gap between Algorithm \ref{DP-DDA} and \cite{LIU202243} is more significant for the case with smaller $\varepsilon$, i.e., a tighter DP requirement.




\subsection{Results for $l_1$-regularized SVM}
Set
$
h(x) = \phi\lVert x \rVert_1
$
with $\phi = 0.0005$. In this case, the problem in \eqref{eq:svm_l1} is convex with a non-smooth regularization term. According to Corollary \ref{C_rate}, we set $\gamma_t= 0.01\sqrt{t}$ and $a_t=1$ in the experiment. 



\begin{figure}[htbp]
\centering
\includegraphics[width=3.3in]{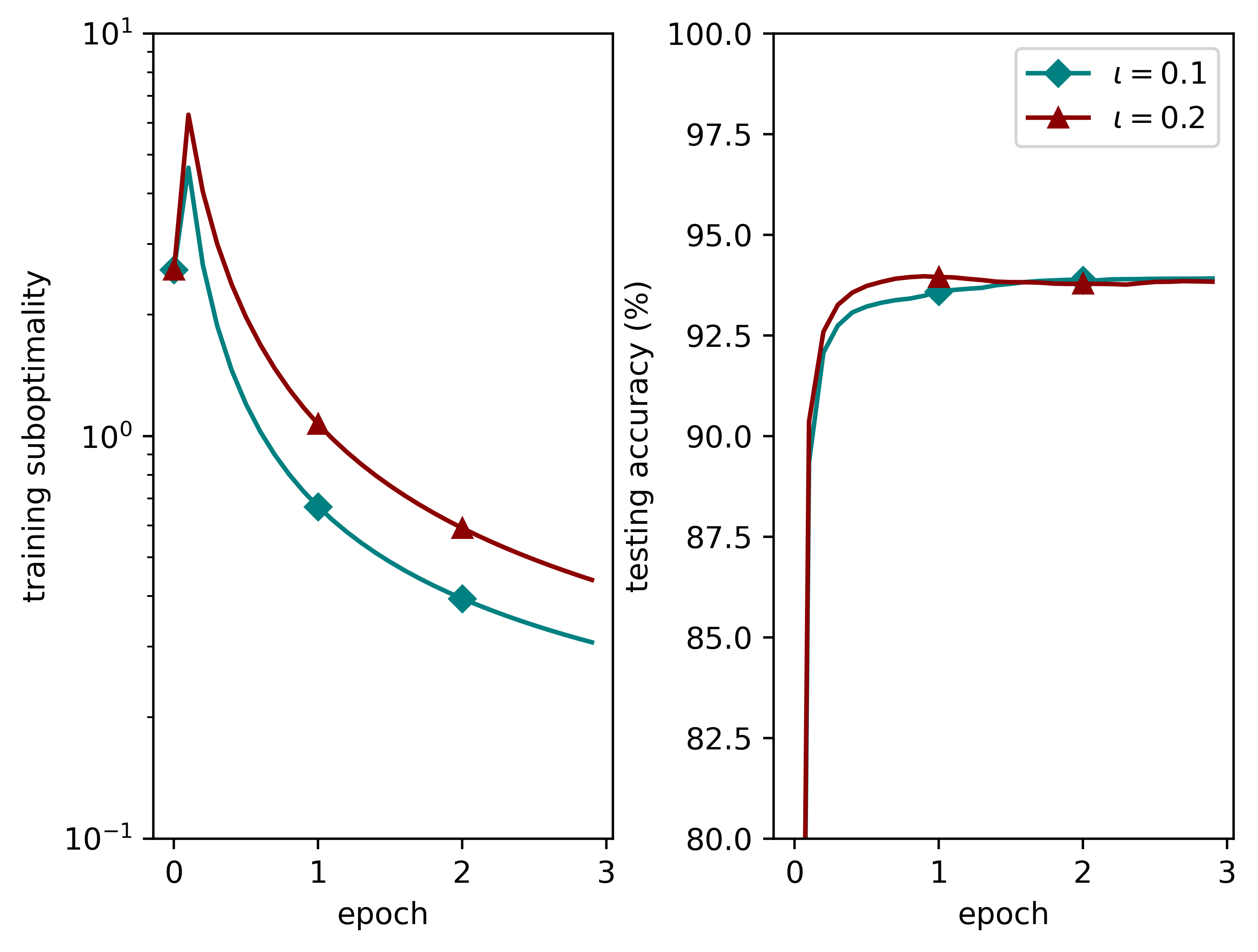}
\caption{
Performance comparison between Algorithm \ref{DP-DDA} with different $\iota$ for $l_1$-regularized SVM with $\varepsilon=0.4$.}
\label{l1_subopt}
\end{figure}

\begin{figure}[htbp]
\centering
\includegraphics[width=3.3in]{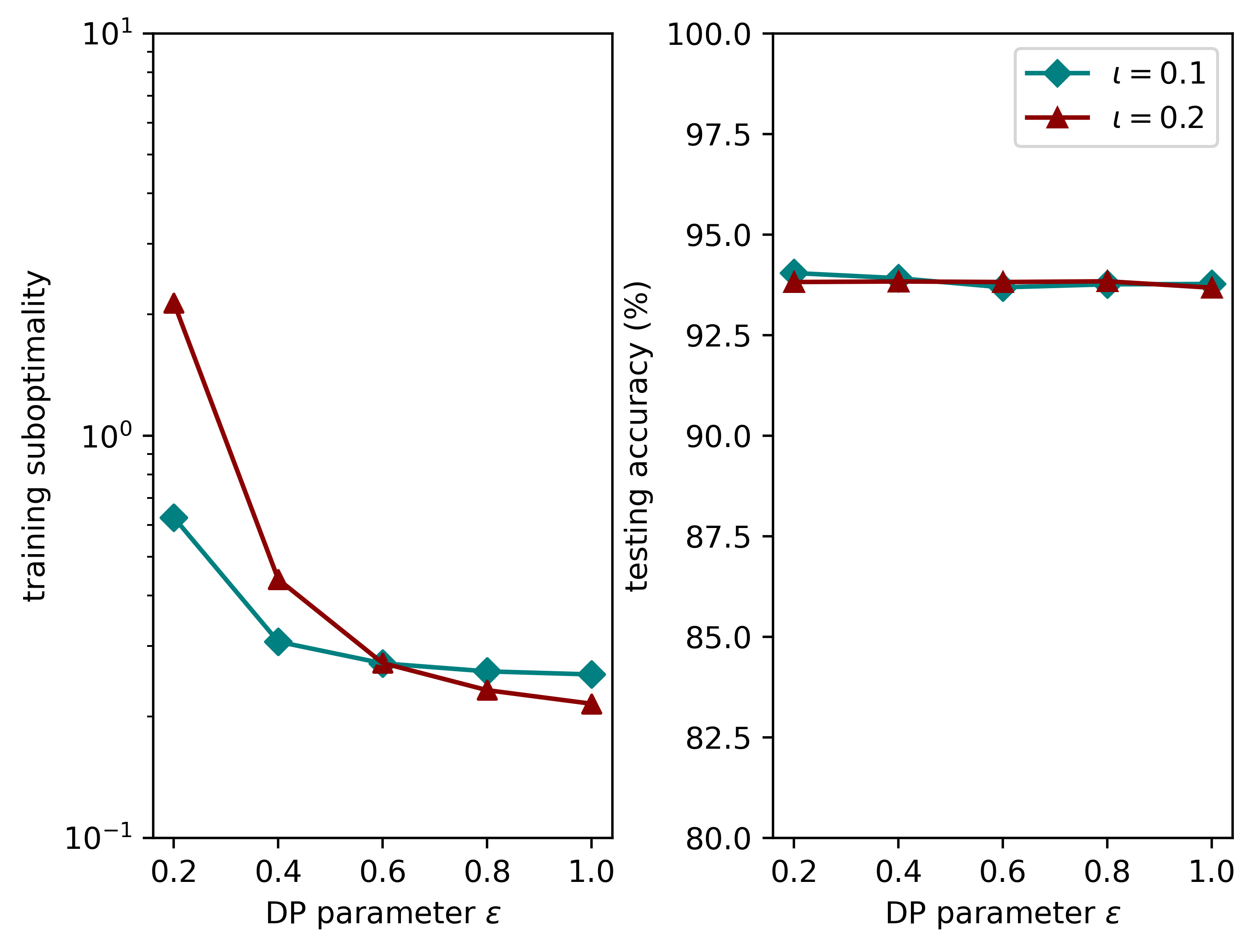}
\caption{Privacy--utility tradeoff in $l_1$-regularized SVM. The suboptimality and accuracy are evaluated after $3$-epoch training.}
\label{l1_nodes_subopt}
\end{figure}

First, we set $\varepsilon=0.4$ and compare  Algorithm \ref{DP-DDA} under different subsampling ratios.
The findings depicted in Fig. \ref{l1_subopt} illustrate a similar trend: As the subsampling ratio $\iota$ decreases, the utility loss diminishes correspondingly.
Additionally, we present the results for Algorithm \ref{DP-DDA} with various DP parameters in Fig. \ref{l1_nodes_subopt}. Notably, for both selected subsampling ratios, we observe a degradation in utility as the DP parameter $\varepsilon$ decreases.



In summary, the experimental results reveal the effectiveness of the proposed algorithms and validate our theoretical findings.

\section{Conclusion}\label{Sec:Concl}
In this work, we presented a class of differentially private DDA algorithms for solving ERM over networks. The proposed algorithms achieve DP by \emph{i}) randomly activating a fraction of nodes at each time instant and \emph{ii}) perturbing the stochastic subgradients over individual data samples within activated nodes. We proved that our algorithms substantially improve over existing ones in terms of utility loss.

{There are numerous promising directions for future endeavors. Firstly, an intriguing avenue to explore is the heterogeneous case, where nodes exhibit substantial variations in dataset size and/or Lipschitz constants. Secondly, it is worthwhile to investigate the high probability convergence of the proposed algorithms.}

\bibliographystyle{plain}        
\bibliography{autosam.bib}           



\appendix

\section{Proof of Lemma \ref{consensus_error}}
\label{Appen:Lem4}

{\bf Notation}: To facilitate the presentation, we introduce the following notation. Define ${\bm W}^{(t)} = W^{(t)}\otimes I$.
Given a real-valued random vector $x$, we let
\begin{equation}\label{expectation_operator}
	\lVert x \rVert_{\mathbb{E}} = \sqrt{  \mathbb{E} [\lVert x \rVert^2]}.
\end{equation}
Accordingly, for a square random matrix $W$, we denote
$
\lVert W \rVert_{\mathbb{E}} =\sup_{\lVert x \rVert_{\mathbb{E}}=1}\lVert Wx \rVert_{\mathbb{E}}
$. 
Denote ${\bf y}^{(t)}	= {\bf  1}\otimes y^{(t)}$ and 
\begin{equation}\label{additional_notation}
	{\bf z}^{(t)}	=\begin{bmatrix}
			z_1^{(t)} \\ \vdots \\	z_n^{(t)}
		\end{bmatrix},		{\bf x}^{(t)}	=\begin{bmatrix}
			x_1^{(t)} \\ \vdots \\	x_n^{(t)}
		\end{bmatrix}, 	
			{{\bm \zeta}}^{(t)}	=\begin{bmatrix}
			{ \zeta}_1^{(t)} \\ \vdots \\	{\zeta}_n^{(t)}
		\end{bmatrix},
			{\bm \eta}^{(t)}	=\begin{bmatrix}
			{ \eta}_1^{(t)} \\ \vdots \\	{\eta}_n^{(t)}
		\end{bmatrix}.
	\end{equation}

\begin{pf}
This proof consists of three parts. First, we prove	
\begin{equation} \label{XYrelation}
		\begin{split}
	\lVert {\bf x}^{(t)} - {\bf y}^{(t)} \rVert  \leq  \frac{1}{\gamma_t+\mu \iota A_t} \left\lVert {\bf z}^{(t)}-{\bm1}\otimes \overline{z}^{(t)}\right\rVert,
		\end{split}
	\end{equation} 
where $\overline{z}(t)=n^{-1}\sum_{i=1}^{n}z_i(t)$.
Second, we prove
\begin{equation}\label{MSE_bound_z}
	  \left\lVert  \frac{{\bf z}^{(t)}}{a_t}-\frac{{\bf 1}\otimes\overline{z}^{(t)}}{a_t} \right\rVert^2_{\mathbb{E}} \leq\frac{n\iota(L^2+m\sigma^2)}{(1-\beta)^2}.
\end{equation}
Finally, we conclude the proof using these two inequalities.

{\bf Part i}) Following \cite[Lemma 5]{liu2021decentralized}, we have $\iota A_\tau h(x)+ \gamma_{\tau} d(x)$ is $(\mu \iota A_\tau+\gamma_\tau)$-strongly convex, and $\Psi_{\tau}^*$, defined by 
\begin{equation}\label{Fenchel_conjugate}
	\begin{split}
	\Psi_{\tau}^* (w) =  \sup_{{x}} \left\{  \langle w, x \rangle - \iota A_\tau h(x)- \gamma_\tau d(x) \right\},
	\end{split}
\end{equation}
has $(\mu \iota A_\tau+\gamma_{\tau})^{-1}$-Lipschitz continuous gradients. Let ${\theta}^{(t)} = n^{-1}\sum_{i=1}^n \eta_i^{(t)}{\zeta}_i^{(t)}$ and $\circ$ be the Hadamard product.
Due to 
\begin{equation*}
	\begin{split}
		&\overline{z}^{(t+1)} = \frac{1}{n}({\bf 1}^T\otimes I) {\bf z}^{(t+1)} \\
		& = \frac{{\bf 1}^T\otimes I}{n}{\bm W}^{(t)} \left({\bf z}^{(t)} + a_t {\bm \eta}^{(t)} \circ {\bm{ \zeta}}^{(t)}  \right)  \\
		& = \frac{({\bf 1}^T\otimes I)}{n} \left({\bf z}^{(t)} + {a_t}{\bm \eta}^{(t)} \circ{\bm{ \zeta}}^{(t)} \right)  \\
		& = \overline{z}^{(t)}+ a_t \theta^{ (t)},
	\end{split}
\end{equation*}
we have
$y^{(\tau)} = \nabla \Psi_{\tau}^*\left(-\sum_{k=1}^{\tau-1} a_k \theta^{(k)}\right) = \nabla \Psi_{\tau}^*\left(-z^{(t)}\right).$
In addition, 
$  x_i^{(\tau)} = \nabla \Psi_{\tau}^*(-{z}_i^{(\tau)}), \, \forall i =1,\dots, n, $
which gives us \eqref{XYrelation}.


%
{\bf Part ii}) When $t=1$, since $z_i^{(1)}=0$ for all $i$, we have $\overline{z}^{(1)}= 0$ and therefore \eqref{MSE_bound_z} satisfied. Next, we consider the case with $t\geq 1$.

Let $\tilde{\bf z}^{(t+1)}= {\bf z}^{(t+1)}-{\bf 1}\otimes \overline{z}^{(t+1)}$.
Then,
\begin{equation}\label{recursion_z}
	\begin{split}
	{\tilde{\bf z}^{(\tau+1)}} =    \tilde{\bm W}^{(\tau)}\left( {\tilde{\bf z}^{(\tau)}}  +  a_\tau\left( {\bm \eta}^{(\tau)} \circ {\bm{ \zeta}}^{(\tau)} -{\bf 1}\otimes \theta^{(\tau)} \right) \right),
\end{split}
\end{equation} 
	where $\tilde{\bm W}^{(\tau-1)}= \tilde{W}^{(\tau-1)}\otimes I$ and $\tilde{W}^{(\tau-1)}= {W}^{(\tau-1)}- \frac{{\bf 1} {\bf 1}^T}{n}$.
By iterating \eqref{recursion_z}, we obtain
\begin{equation}\label{iteration_z}
	\begin{split}
&{\tilde{\bf z}^{(t)}}= \sum_{\tau=1}^{t-1}    \tilde{\bm W}^{(\tau,t-1)}   a_\tau \left( {\bm \eta}^{(\tau)} \circ {\bm{ \zeta}}^{(\tau)} - {\bf 1}\otimes   \theta^{(\tau)} \right)\\
& =  \sum_{\tau=1}^{t-1}    \tilde{\bm W}^{(\tau,t-1)}   a_\tau\left( {\bm \eta}^{(\tau)} \circ {\bm{ \zeta}}^{(\tau)} \right),
\end{split}
\end{equation}
where $\tilde{\bm W}^{(\tau,t-1)}= \tilde{\bm W}^{(t-1)}\dots \tilde{\bm W}^{(\tau)}$.
Since $ {\bf 1}\otimes \theta^{(\tau)} =({\bf 1}\otimes I)\theta^{(\tau)} $ and $(A\otimes B)(C\otimes D) = (AC) \otimes (BD)$, we have
$\left( \tilde{W}^{(\tau,t-1)}\otimes I  \right) ({\bf 1}\otimes \theta^{(\tau)} )  = \left(\left( \tilde{W}^{(\tau,t-1)} {\bf 1}  \right) \otimes I  \right)\theta^{(\tau)}  = 0$. Therefore,
\eqref{iteration_z} can be rewritten as,  $\forall t \geq 1$,
\begin{equation*}
	\begin{split}
	{\tilde{\bf z}^{(t)}} = \sum_{\tau=1}^{t-1}   \tilde{\bm W}^{(\tau,t-1)}   a_\tau \left( {\bm \eta}^{(\tau)} \circ {\bm{ \zeta}}^{(\tau)} \right).
	\end{split}
\end{equation*}
Upon taking the norm defined in \eqref{expectation_operator} on both sides and using the Minkowski inequality~\cite{gut2013probability}, we obtain
\begin{equation}\label{norm_E}
\begin{split}
    & \lVert \tilde{\bf{z}}^{(t)} \rVert_\mathbb{E}\leq  \sum_{\tau=1}^{t-1}  a_\tau  \left \lVert   \tilde{\bm W}^{(\tau,t-1)} \left( {\bm \eta}^{(\tau)} \circ {\bm{ \zeta}}^{(\tau)}\right)\right \rVert_{\mathbb{E}}.
\end{split}
\end{equation}
Consider
\begin{equation*}
\begin{split}
     &\left \lVert  \tilde{\bm W}^{(\tau,t-1)}\left( {\bm \eta}^{(\tau)} \circ {\bm{ \zeta}}^{(\tau)}\right) \right \rVert^2_{\mathbb{E}} \\
     &\overset{(i)}{=} \mathbb{E}\Big[ \left( {\bm \eta}^{(\tau)} \circ {\bm{ \zeta}}^{(\tau)}\right)^T (\tilde{\bm W}^{(\tau,t-2)})^T  \Big((\tilde{\bm W}^{(t-1)})^T  \tilde{\bm W}^{(t-1)}\Big) \\
     & \quad \quad \,\, \times \tilde{\bm W}^{(\tau,t-2)} \left( {\bm \eta}^{(\tau)} \circ{\bm{ \zeta}}^{(\tau)}\right)\Big] \\
     & \overset{(ii)}{\leq} \left\lVert  \tilde{\bm W}^{(\tau,t-2)} \left( {\bm \eta}^{(\tau)} \circ {\bm{ \zeta}}^{(\tau)}\right) \right\rVert^2_{\mathbb{E}} \rho\left(\mathbb{E}_W\left[(\tilde{W}^{(t-1)})^T  \tilde{W}^{(t-1)}   \right]\right) \\
    & \overset{(iii)}{\leq } \beta^{2}\left\lVert    \tilde{\bm W}^{(\tau,t-2)} \left( {\bm \eta}^{(\tau)} \circ {\bm{ \zeta}}^{(\tau)}\right)  \right\rVert^2_{\mathbb{E}}\\
     & \overset{(iv)}{\leq } \beta^{2(t-\tau-1)} \left\lVert \tilde{\bm W}^{(\tau)} \left( {\bm \eta}^{(\tau)} \circ{\bm{ \zeta}}^{(\tau)}\right)  \right\rVert^2_{\mathbb{E}}, 
     \end{split}
\end{equation*}
where we use $(A\otimes B)(C\otimes D)=(AC)\otimes(BD)$ and that $\tilde{W}^{(t-1)}$ is independent of the random events that occur up to time $t-2$ to obtain (i) and (ii), respectively, (iii) is due to
$\mathbb{E}[(\tilde{W}^{(t-1)})^T\tilde{W}^{(t-1)}] = \mathbb{E}[(W^{(t-1)})^T W^{(t-1)}]- \frac{{\bf 1} {\bf 1}^T}{n}$ and (iv) is by iteration. Therefore, we get from \eqref{norm_E} that
\begin{equation*}
\begin{split}
    &\lVert \tilde{\bf{z}}^{(t)} \rVert_\mathbb{E} \leq  \sum_{\tau=1}^{t-1} \beta^{t-\tau-1}  a_\tau  \left\lVert \tilde{\bm W}^{(\tau)} \left( {\bm \eta}^{(\tau)} \circ {\bm{ \zeta}}^{(\tau)}\right)  \right\rVert_{\mathbb{E}}  \\
    & \leq  \sum_{\tau=1}^{t-1} \beta^{t-\tau-1}  a_\tau  \left\lVert \tilde{\bm W}^{(\tau)} \right\rVert_{\mathbb{E}}  \left\lVert {\bm \eta}^{(\tau)} \circ {\bm{ \zeta}}^{(\tau)} \right\rVert_{\mathbb{E}} \\
    &  \leq \sum_{\tau=1}^{t-1} \beta^{t-\tau-1}  a_\tau   \left\lVert {\bm \eta}^{(\tau)} \circ {\bm{ \zeta}}^{(\tau)} \right\rVert_{\mathbb{E}},
\end{split}
\end{equation*}
where the last inequality is due to $\left\lVert \tilde{\bm W}^{(\tau)} \right\rVert_{\mathbb{E}}  \leq 1$,
leading to
\begin{equation*}
\begin{split}
    &\lVert \tilde{\bf{z}}^{(t)} \rVert^2_{\mathbb{E}} \leq \left(\sum_{\tau=1}^{t-1}  a_\tau \beta^{t-\tau-1} \left\lVert {\bm \eta}^{(\tau)} \circ {\bm{ \zeta}}^{(\tau)} \right\rVert_{\mathbb{E}} \right)^2 \\
    & \leq  \left(\sum_{\tau=1}^{t-1}  \left(\beta^{\frac{t-\tau-1}{2}} \right)^2 \right) \left(\sum_{\tau=1}^{t-1} \left( a_\tau \beta^{\frac{t-\tau}{2}} \left\lVert {\bm \eta}^{(\tau)} \circ {\bm{ \zeta}}^{(\tau)} \right\rVert_{\mathbb{E}} \right)^2\right) \\ 
    & \leq  \frac{1}{1-\beta} \sum_{\tau=1}^t a_\tau^2\beta^{t-\tau-1} \left\lVert {\bm \eta}^{(\tau)} \circ {\bm{ \zeta}}^{(\tau)} \right\rVert^2_{\mathbb{E}},
\end{split}
\end{equation*}
where the second inequality is due to $(w+v)^2\leq 2w^2+2v^2$. Upon dividing both sides by $a_t^2$ and using $0<a_{\tau}\leq a_t, \forall \tau\leq t$, we have
\begin{equation*}
\begin{split}
    &\left\lVert \frac{\tilde{\bf{z}}^{(t)}}{a_t} \right\rVert^2_{\mathbb{E}}  \leq  \frac{1}{1-\beta} \sum_{\tau=1}^t \frac{a_\tau^2}{a_t^2}\beta^{t-\tau-1} \left \lVert {\bm \eta}^{(\tau)} \circ {\bm{ \zeta}}^{(\tau)}   \right \rVert^2_{\mathbb{E}} \\
    &\leq \frac{n \iota (L^2+m \sigma^2)}{(1-\beta)^2}.
\end{split}
\end{equation*}


\noindent{\bf Part iii}) 
Based on~\eqref{XYrelation} and~\eqref{MSE_bound_z},
we have
\begin{equation*}
\begin{split}
  \frac{{n}\iota(L^2+m\sigma^2)a_t^2}{(1-\beta)^2(\mu \iota A_t+\gamma_t)^2}\geq \lVert {\bf x}^{(t)}- {\bf y}^{(t)} \rVert^2_{\mathbb{E}} 
  \geq  \sum_{i=1}^n \lVert {x}_i^{(t)}- {y}^{(t)} \rVert^2_{\mathbb{E}} .
  \end{split}
\end{equation*}
By the Jensen's inequality, we have
\begin{equation*}
\begin{split}
     \mathbb{E}\left[ \lVert  {\bf x}^{(t)}-{\bf y}^{(t)} \rVert \right]   \leq \sqrt{\lVert  {\bf x}^{(t)}-{\bf y}^{(t)} \rVert^2_{\mathbb{E}}
     }\leq \frac{\sqrt{n\iota}(L+\sqrt{m}\sigma)a_t}{(\gamma_t+\mu \iota A_t)(1-\beta)}.
\end{split}
\end{equation*}
This together with the bound between $l_1$ and $l_2$-norms, i.e.,
$
 \frac{1}{ \sqrt{n}} \sum_{i=1}^n\lVert x_i^{(t)}-y^{(t)} \rVert  \leq  \lVert  {\bf x}^{(t)}-{\bf y}^{(t)} \rVert 
$
yields
\begin{equation}\label{bound_z}
\sum_{i=1}^n \mathbb{E}[\lVert x_i^{(t)}-y^{(t)} \rVert] \leq 	\frac{{n\iota}(L+\sqrt{m}\sigma)a_t}{(\gamma_t+\mu \iota A_t)(1-\beta)}.
\end{equation}
This completes the proof.
\end{pf}



\section{Proof of Lemma \ref{auxiliary_sequence}}
\label{Appen:Lem5}
\begin{pf}
Recall \eqref{Fenchel_conjugate}
\begin{equation*}
	\begin{split}
		\Psi_{\tau}^* (w) =  \sup_{{x}} \left\{  \langle w, x \rangle - \iota A_\tau h(x)- \gamma_{\tau} d(x) :=R_{\tau}(x,w) \right\}.
	\end{split}
\end{equation*}
	Since $-R_{\tau}(\cdot,w)$ is $(\mu \iota  A_\tau+\gamma_\tau)$-strongly convex, we have $$\nabla \Psi_{\tau}^*\left(-  \sum_{k=1}^{\tau-1} a_k {\theta}^{(k)}\right)=y^{(\tau)},$$ where  ${\theta}^{(t)} = n^{-1}\sum_{i=1}^n \eta_i^{(t)}{\zeta}_i^{(t)}$,  and $\Psi_{{\tau}}^*$ has $(\mu \iota A_\tau+\gamma_{\tau})^{-1}$-Lipschitz continuous gradients, see, e.g., \cite[Lemma 5]{liu2021decentralized},
implying 
	\begin{equation}\label{DA-smoothness}
	\begin{split}
		&\Psi_{{\tau}}^*\left(-  \sum_{k=1}^{\tau} a_k {\theta}^{(k)}\right) \\
		&\leq 		\Psi_{{\tau}}^*\left(- \sum_{k=1}^{\tau-1} a_k {\theta}^{(k)}\right)   -{a_{\tau}}\langle {y}^{(\tau)} ,  {\theta}^{(\tau)} \rangle
  +\frac{a_\tau^2 \left\lVert   {\theta}^{(\tau)} \right\rVert^2}{2(\mu \iota A_\tau+\gamma_{\tau})} .
	\end{split}
\end{equation}
Upon using $\gamma_{\tau}$ is non-decreasing and $d(x)\geq 0$, we have
	\begin{equation}\label{DA-optimality}
			\begin{split}
				&\Psi_{{\tau}}^*\left(- \sum_{k=1}^{\tau-1} a_k {\theta}^{(k)}\right)  =R_{\tau}\left(y^{(\tau)},- \sum_{k=1}^{\tau-1} a_k  {\theta}^{(k)}\right)\\ &= R_{\tau-1}\left(y^{(\tau)}, - \sum_{k=1}^{\tau-1} a_k  {\theta}^{(k)}\right) -\iota a_{\tau}h(	y^{(\tau)} )\\
				&\quad + ( \gamma_{\tau-1} -\gamma_{\tau})d(	y^{(\tau)} )\\
				& \leq R_{\tau-1}\left(y^{(\tau)}, - \sum_{k=1}^{\tau-1} a_k  {\theta}^{(k)}\right)-\iota a_{\tau}h(	y^{(\tau)} ) \\
				& \leq \Psi_{{\tau-1}}^*\left(-\iota \sum_{k=1}^{\tau-1} a_k  {\theta}^{(k)}\right) -\iota a_\tau h(	y^{(\tau)} ).
			\end{split}
		\end{equation}
Upon plugging \eqref{DA-optimality} into \eqref{DA-smoothness} and summing up the resultant inequality from $\tau=1$ to $\tau=t$, we have
	\begin{equation*}
		\begin{split}
	&	 \sum_{\tau=1}^{t} a_{\tau}	\left(\langle  {\theta}^{(\tau)}, {y}^{(\tau)} \rangle +  \iota h(y^{(\tau)})\right)	\\
	&	\leq  \Psi_{0}^*(0) -\Psi_{\tau}^*\left(- \sum_{k=1}^{\tau} a_k  {\theta}^{(k)} \right)     +\sum_{\tau=1}^{t}\frac{a_{\tau}^2  \lVert  {\theta}^{(\tau)} \rVert^2 }{2(\mu \iota A_\tau+\gamma_{\tau})}.
				\end{split}
	\end{equation*}
	Note that $y^{(1)}= \nabla \Psi_{{1}}^*(0)$, $A_0=0$ and $\gamma_0=0$ by definition, implying that $\Psi_{0}^*(0) =0$.
	Further considering
	\begin{equation*}
		\begin{split}
		& \sum_{\tau=1}^{t} a_{\tau}\langle\theta ^{(\tau)}, -x^*\rangle \leq \iota A_t h(x^*)+\gamma_{t} d(x^*)\\
		& \quad + \sup_{ x  } \left\{  	- \sum_{k=1}^{\tau} a_k \langle \theta^{(k)}, x\rangle - \iota A_{t} h(x)-\gamma_{t}d(x)\right\}  \\
		&\leq 	\Psi_{t}^* \left(- \sum_{k=1}^{\tau} a_k \theta^{(k)}\right) + \iota A_t h(x^*)+ \gamma_{t} d(x^*),
		\end{split}
	\end{equation*}
	we obtain
	\begin{equation}
		\begin{split}
		& \sum_{\tau=1}^{t} a_\tau	\Big( \langle \theta^{(\tau)}, {y}^{(\tau)}-x^* \rangle   + \iota\left(h(y^{(\tau)}) -h(x^*) \right) \Big)	 \\
		&\leq   \sum_{\tau=0}^{t}\frac{a_\tau^2 \lVert \theta^{(\tau)} \rVert^2  }{2(\mu \iota A_\tau+\gamma_{\tau})} + \gamma_{t} d(x^*).
	\end{split}
	\end{equation}
	Dividing both sides by $\iota>0$ leads to the desired inequality.
\end{pf}

\section{Proof of Corollary \ref{SC_rate}}
\label{Cor:SC_rate}

\begin{pf}
We obtain from the update of $A_t$ in Algorithm \ref{DP-DDA} that
\begin{equation}\label{SC_scheduling}
	\sum_{\tau=1}^{t} \frac{a_\tau^2}{\mu \iota A_\tau+\gamma_\tau} = \sum_{\tau=1}^{t} \frac{2\tau^2}{\mu \iota \tau(\tau+1)} \leq  \frac{2t}{\mu\iota }.
\end{equation}
Due to $\mu$-strong convexity of $F$, we have
\begin{equation}\label{x_to_xstar}
	\begin{split}
		& \sum_{i=1}^n\lVert  \tilde{x}_i^{(t)} -  x^*  \rVert^2 \leq  	2 \sum_{i=1}^n\lVert  \tilde{x}_i^{(t)} -  \tilde{y}^{(t)}  \rVert^2  + 2n	\lVert    \tilde{y}^{(t)} - x^* \rVert^2  \\
		&\leq   	2\lVert  \tilde{x}_i^{(t)} -  \tilde{y}^{(t)}  \rVert^2  + \frac{4n}{\mu } \Big( F(\tilde{y}^{(t)})-F(x^*)\Big).
	\end{split}
\end{equation}
Upon using convexity of the $l_2$-norm and \eqref{MSE_consensus_expectation}, we obtain
\begin{equation}\label{MSE_consensus}
	\begin{split}
		& \sum_{i=1}^n \mathbb{E}\left[\lVert  \tilde{x}_i^{(t)} -  \tilde{y}^{(t)} \rVert^2 \right]\leq \mathbb{E}\left[\frac{1}{A_t }\sum_{\tau=0}^{t} \sum_{i=1}^n a_\tau	\lVert   x_i^{(\tau)}-y^{(\tau)}\rVert^2 \right] \\
		& \leq  \frac{1}{A_t } \sum_{\tau=1}^{t}  \frac{a_\tau^3}{(\mu \iota A_\tau+\gamma_\tau)^2}\left(  \frac{\iota nL^2}{(1-\beta)^2} \right) \\
		& \leq \frac{4nL^2}{\mu^2 \iota A_t(1-\beta)^2}  \sum_{\tau=1}^{t} \frac{1}{a_\tau} \leq \frac{4nL^2 (\log t+1)}{\mu^2 \iota A_t(1-\beta)^2}  .
	\end{split}
\end{equation}
By putting expectations on \eqref{x_to_xstar} and using \eqref{MSE_consensus} and \eqref{error_bound}, we arrive at \eqref{rate_mu_convex} as desired.
\end{pf}

\section{Proof of Corollary \ref{C_rate}}
\label{Cor:C_rate}

\begin{pf}
	Under the choices of $a_t=1$ and $\gamma_t=\gamma\sqrt{t}$, we have
	\begin{equation*}
		\sum_{\tau=1}^{t} \frac{a_\tau^2}{\mu \iota A_\tau +\gamma_t} =\frac{1}{\gamma} \sum_{\tau=1}^{t} \frac{1}{\sqrt{\tau}}\leq \frac{2}{\gamma}\sqrt{t}.
	\end{equation*} 
	Therefore, by using \eqref{error_bound} and \eqref{error_consensus}, we obtain
	for all $t\geq 1$,
	\begin{equation*}
		\mathbb{E}\left[ F(\tilde{y}^{(t)})-F(x^*) \right]\leq\frac{d(x^*) +2\iota M}{\iota \gamma\sqrt{t}} 
	\end{equation*}
	and
	\begin{equation*}
		\frac{1}{n}\sum_{i=1}^n\mathbb{E}\left[	\lVert   x_i^{(t)}-y^{(t)}\rVert \right] \leq \frac{L\sqrt{\iota}}{\gamma(1-\beta)\sqrt{t}}, \quad i=1,\dots, n
	\end{equation*}
	respectively. Since $l_2$-norm is convex, we have
	\begin{equation}\label{Eq:ergodic_average}
		\frac{1}{n}\sum_{i=1}^n\mathbb{E}[\lVert \tilde{x}_i^{(t)} - \tilde{y}^{(t)}\rVert] \leq \frac{1}{t}\sum_{\tau=1}^t\frac{L\sqrt{\iota} }{\gamma(1-\beta)\sqrt{t}} \leq  \frac{2L\sqrt{\iota}}{\gamma(1-\beta)\sqrt{t}}.
	\end{equation}
	
\end{pf}

\section{Proof of Corollary \ref{Cor:SC_loss}}
\label{App:SC_loss}

\begin{pf}
	By \eqref{SC_scheduling}, we have
		\begin{equation*}
		\begin{split}
		&	\mathbb{E}[F(\tilde{y}^{(T)})-F(x^*)] \leq \frac{4M}{\mu\iota(T+1)}+ \frac{80m \log(2/\delta_0)\iota^2 L^2 }{\mu  q^2\varepsilon^2} \\
		& \quad \quad \quad \quad \quad \quad \quad \quad \quad\quad +\frac{16\sqrt{10m\log(2/\delta_0)\iota}L^2}{\mu  q\varepsilon (1-\beta)\sqrt{T+1}}  ,  
		\end{split}
	\end{equation*}
	where $M$ is a positive constant defined in Theorem \ref{thm_error_bound}.
	Similar to \eqref{MSE_consensus}, we have
	\begin{equation*}
	\begin{split}
		 \frac{1}{n}\sum_{i=1}^n \mathbb{E}\left[\lVert  \tilde{x}_i^{(T)} -  \tilde{y}^{(T)} \rVert^2 \right] \leq \frac{8(L^2+m\sigma^2) (\log T+1)}{\mu^2 \iota T(T+1)(1-\beta)^2}  .
	\end{split}
\end{equation*}
Following \eqref{x_to_xstar}, we set $T = \mathcal{O}\left(\frac{q^2\varepsilon^2}{\iota^3  (1-\beta)^2m\log(1/\delta_0)}\right)$ to obtain \eqref{ULoss:SC}.
\end{pf}

\section{Proof of Corollary \ref{Cor:C_loss}}
\label{App:C_loss}

\begin{pf}
	By the specific choices of parameters in \eqref{parameter_choices_C_DP}, we obtain from \eqref{error_bound} that
		\begin{equation*}
		\begin{split}
			&\mathbb{E}[F(\tilde{y}^{(T)})-F(x^*)] \\
			&\leq  \frac{d(x^*) +2\iota M}{\gamma \iota \sqrt{T}}+ \frac{40m\iota^3 L^2 \log(2/\delta_0)}{\gamma {q^2}\varepsilon^2 } \sum_{t=1}^{T} \frac{1}{\sqrt{t}}  \\
			& \quad + \frac{4\sqrt{10 m \iota \log(2/\delta_0)}\iota L^2}{\gamma q\varepsilon(1-\beta) \sqrt{T}}  \sum_{t=1}^{T} \frac{1}{\sqrt{t}}  \\
			& \leq \frac{d(x^*) +2\iota M}{\gamma \iota \sqrt{T}} + \frac{8\sqrt{10 m \iota \log(2/\delta_0)}\iota L^2}{\gamma q\varepsilon(1-\beta) }  \\
		&	\quad + \frac{80  m \iota^3 L^2 \log(2/\delta_0) \sqrt{T}}{ {\gamma}{q^2}\varepsilon^2},
		\end{split}
	\end{equation*}
where the last inequality is due to
$
	\sum_{t=1}^{T} \frac{1}{\sqrt{t}} \leq 1+ \int_{1}^{T}\frac{1}{\sqrt{t}} dt \leq 2\sqrt{t}.
$
Upon using $  \iota \leq 1-\beta$ and $T=\mathcal{O}\left(\frac{q^2\varepsilon^2}{\iota^3(1-\beta)^2{m\log(1/\delta_0)}}\right)$, we arrive at \eqref{ULoss:C}. In addition, similar to \eqref{Eq:ergodic_average}, we have
\begin{equation*}
\begin{split}
	&	\frac{1}{n}\sum_{i=1}^n\mathbb{E}[\lVert \tilde{x}_i^{(T)} - \tilde{y}^{(T)}\rVert] \\
 & \leq  \frac{2(L+\sqrt{m}\sigma)\sqrt{\iota}}{\gamma(1-\beta)\sqrt{T}} \leq  \mathcal{O}\left(\frac{L\sqrt{m\iota \log(1/\delta_0)}}{\gamma q\varepsilon}\right).
\end{split}
\end{equation*}
\end{pf}

\end{document}